\documentclass[12pt,a4paper,draft]{article}
\usepackage[T2A]{fontenc}
\usepackage[utf8]{inputenc}
\usepackage[english,russian]{babel}
\usepackage{srcltx}
\usepackage{cite}
\usepackage{latexsym,amsfonts,amssymb,amsmath,longtable,stmaryrd,mathrsfs,amsthm,hyperref}
\usepackage{stmaryrd} 
\usepackage{calrsfs}
\usepackage[mathlines,pagewise]{lineno}\usepackage{multirow}
\renewcommand{\leq}{\leqslant}
\renewcommand{\geq}{\geqslant}

\newcommand{\la}{\langle}
\newcommand{\ra}{\rangle}
\newcommand{\BS}{{
\mathcal B}{\mathcal S}}
\theoremstyle{plain}
\newtheorem{Lemma}{{\bfseries Лемма}}

\newtheorem{Theo}{{\bfseries Теорема}}

\DeclareMathOperator{\Irr}{Irr} 
\DeclareMathOperator{\GL}{GL} \DeclareMathOperator{\Aut}{Aut} \DeclareMathOperator{\Inn}{Inn}
\DeclareMathOperator{\SL}{SL} 
 \DeclareMathOperator{\Sp}{Sp}
\DeclareMathOperator{\GU}{GU} 
\DeclareMathOperator{\Oo}{O} \DeclareMathOperator{\SU}{SU}
\DeclareMathOperator{\PSL}{PSL} \DeclareMathOperator{\PSU}{PSU}
\newcommand{\F}{\mathbb F}

\renewcommand{\L}{\widehat{L}}

\DeclareMathOperator{\Sym}{Sym}
\DeclareMathOperator{\bs}{BS}

\title{\vspace{-1cm} 
On the sharp Baer--Suzuki theorem for the $\pi$-radical\thanks{Revin  is supported by RFBR and BRFBR, project \textnumero\ 20-51-00007. Revin and Vdovin are supported by the
State Contract of the Sobolev Institute of Mathematics, project \textnumero\ 0314-2019-0001. }}
\begin{document}

\author{ Nanying Yang, Zhenfeng Wu\\
{\small Jiangnan University,}\\ {\small Wuxi 230026, P. R. China}\\
{\small E-mail:
yangny@jiangnan.edu.cn, zhfwu@jiangnan.edu.cn}\\ \\
Danila O. Revin, Evgeny P. Vdovin\\
{\small Sobolev Institute of Mathematics and Novosibirsk State University,}\\
{\small Novosibirsk 630090, Russia}\\
{\small E-mail: revin@math.nsc.ru, vdovin@math.nsc.ru}}

 \date{}
\maketitle

\pagenumbering{arabic}

{\small 

\centerline{\bf Abstract}

\medskip

Let $\pi$ be a set of primes such that $|\pi|\geqslant 2$ and $\pi$ differs from the set of all primes. Denote by $r$ the smallest prime which does not belong to $\pi$ and set $m=r$ if $r=2,3$ and $m=r-1$ if $r\geqslant 5$. We study the following conjecture: a conjugacy class $D$ of a finite group $G$ is contained in $\Oo_\pi(G)$ if and only if every $m$ elements of $D$ generate a $\pi$-subgroup. We confirm this conjecture for each group $G$ whose nonabelian composition factors are isomorphic to alternating, linear and unitary simple groups.


\medskip
\noindent {\bf Key words:} linear and unitary finite simple groups, $\pi$-radical, Baer--Suzuki $\pi$-theorem.

\medskip
\noindent {\bf MSC2010:} 20D20
}

\section*{Введение}

В работе всюду через $\pi$ обозначается некоторое множество простых чисел.
Конечная группа называется {\it $\pi$-группой}, если все простые делители ее порядка принадлежат~$\pi$. Используются следующие стандартные обозначения. Для конечной группы $G$ через $\Oo_\pi(G)$ обозначается {\it $\pi$-радикал}, т.~е. наибольшая нормальная $\pi$-подгруппа группы $G$. Если $M$~--- подмножество группы $G$, то через $\langle M\rangle$ обозначена подгруппа, порождённая множеством $M$.
Для любой группы $G$ через $G^\sharp$ обозначено множество $G\setminus\{1\}$.

Теорема Бэра-Судзуки \cite{Baer,Suz,AlpLy} утверждает

\medskip\noindent
{\bf Теорема Бэра-Судзуки.}  {\it Пусть $p$ --- некоторое
простое число, $G$ ---  конечная группа и $x\in G$. Тогда
$x\in \Oo_p(G)$ если и только если  $\la x,x^g \ra$ является $p$-группой для
любого $g\in G$. }
\medskip

Нетривиальным в этой теореме является утверждение ``если''. Теорема Бэра--Су\-дзу\-ки имеет несколько эквивалентных формулировок.
Обобщения и аналоги теоремы Бэра-Судзуки исследовались многими авторами, см., напр.,~\cite{AlpLy,MOS,Mamont,Soz,Fl,Gu,FGG,GGKP,GGKP1,GGKP2, Palchik, Tyut,Tyut1,BS_odd,BS_Dpi}. Например, Н.\,Гордеевым, Ф.\,Грюневальдом, Б.\,Кунявским и Е.\,Плоткиным \cite{GGKP2} и независимо П.\,Флавеллом, С.\,Гэстом и Р.\,Гуральником \cite{FGG} доказано, что если любые четыре элемента из данного класса сопряженности в конечной группе порождают разрешимую группу, то  весь класс содержится в разрешимом радикале группы.

В работе \cite[теорема~1]{YRV} авторами было доказано, что для любого множества $\pi$ простых чисел существует число $m=m(\pi)$ такое, что в произвольной группе $G$ выполнено равенство
$$
\Oo_\pi(G)=\{x\in G\mid \langle x^{g_1}, \dots, x^{g_m}\rangle~\text{--- } \pi\text{-группа для любых } g_1,\dots g_m\in G\}.
$$ Наименьшее такое число $m$ согласно \cite{GGKP1} называется шириной Бэра--Сузуки класса $\pi$-групп и обозначается $\bs(\pi)$. Кроме того, доказано \cite[теорема~2]{YRV}, что если $\pi$~--- непустое собственное подмножество множества всех простых чисел и $r$~--- наименьшее простое число, не  лежащее в~$\pi$, то
$$
r-1\leq \bs(\pi)\leq \max\{11, 2(r-2)\}.
$$
Там же было высказано предположение, что нижняя оценка $r-1$ в большинстве случаев совпадает с~$\bs(\pi)$. Более точно:

\medskip\noindent
{\bf Гипотеза 1.} \cite{YRV}.  { Пусть $\pi$ --- собственное подмножество множества всех простых чисел, содержащее как минимум два простых числа, и $r$~--- наименьшее простое число, не входящее в~$\pi$. Тогда $$\bs(\pi)\leq\left\{\begin{array}{rl}
                               r, & \text{ если } r\in\{2,3\}, \\
                               r-1, &   \text{ если } r\geq 5.
                             \end{array}\right.
$$}

В.\,Н.\,Тютянов \cite{Tyut1} подтвердил гипотезу~1 при $r=2$.  Известны также примеры \cite[пример~2]{BS_odd} множеств $\pi$ с $r=3$ таких, что $\bs(\pi)>2$. Результаты \cite{BS_odd} сводят гипотезу к изучению почти простых групп. Чтобы сформулировать утверждение о почти простых группах, которое гарантировало бы справедливость гипотезы 1, напомним обозначения, введенные в~\cite{GS} и~\cite{YRV}.

Пусть $L$~--- неабелева простая группа, $r$~--- простой делитель её порядка, и $x\in\Aut(L)^\sharp$~--- ее автоморфизм. Мы отождествляем $L$ c подгруппой $\Inn(L)$ в~$\Aut(L)$. Cогласно~\cite{GS} через $\alpha(x,L)$ (или просто $\alpha(x)$, если группа $L$ понятна из контекста) обозначено наименьшее количество $L$-сопряжённых с $x$ элементов, которые порождают $\langle L,x\rangle$. По аналогии c $\alpha(x,L)$ авторы в работе~\cite{YRV} обозначили через $\beta_r(x,L)$ (или просто $\beta_r(x)$, если из контекста ясно, что понимается под группой $L$) наименьшее количество $L$-сопряжённых с $x$ элементов, которые порождают подгруппу группы $\langle L,x\rangle$, порядок которой делится\footnote{Если $r$ не делит $|L|$, то величина $\beta_r(x,L)$ не определена.} на~$r$. Ясно, что $\beta_r(x,L)\leqslant \alpha(x,L)$.

Учитывая имеющуюся редукцию к почти простым группам и результаты для $r=2$, см. \cite{Tyut1,BS_odd}, справедливость гипотезы 1 была бы доказана, если бы удалось установить справедливость следующего утверждения:

\medskip\noindent
{\bf Гипотеза 2.} \cite{YRV}.
Пусть $L$~--- неабелева конечная простая группа, $x\in \Aut(L)$~--- элемент простого порядка. Пусть также $r$~--- нечётное простое число, и простой делитель  $s$ порядка группы $L$ выбран так, что ${s=r}$, если $r$ делит $|L|$, и $s> r$ в противном случае. Тогда
  $$
  \beta_{s}(x,L)\leq
  \left\{
  \begin{array}{rl}
  3,&\text{если }r=3,\\
  r-1,&\text{если }r>3.
  \end{array}
  \right.
  $$

Для случая, когда $L$~--- знакопеременная группа, гипотеза 2 верна \cite[предложение~2]{YRV}. В данной статье мы подтвердим гипотезу 2 для простых линейных и унитарных групп. Таким образом, справедлива

\begin{Theo}\label{main}
Пусть $L=L_n(q)$ или $L=U_n(q)$~--- конечная простая линейная или унитарная группа, $x\in \Aut(L)$~--- элемент простого порядка. Пусть также $r$~--- нечётное простое число, и простой делитель  $s$ порядка группы $L$ выбран так, что $s=r$, если $r$ делит $|L|$, и $s> r$ в противном случае.  Тогда
  $$
  \beta_{s}(x,L)\leq
  \left\{
  \begin{array}{rl}
  3,&\text{если }r=3,\\
  r-1,&\text{если }r>3.
  \end{array}
  \right.
  $$
\end{Theo}

С помощью этой теоремы, используя результаты работы \cite{BS_odd}, получаем следующее утверждение.

\begin{Theo}\label{Cor}
Пусть $\pi$~--- некоторое множество простых чисел и $r$~--- наименьшее простое число, не лежащее в $\pi$. Положим
$$m=\left\{\begin{array}{rl}
                               r, & \text{ если } r\in\{2,3\}, \\
                               r-1, &   \text{ если } r\geq 5.
                             \end{array}\right.
$$
Тогда
$$
\Oo_\pi(G)=\{x\in G\mid \langle x^{g_1}, \dots, x^{g_m}\rangle~\text{--- } \pi\text{-группа для любых } g_1,\dots g_m\in G\}
$$
для любой конечной группы $G$, всякий неабелев композиционный фактор которой изоморфен знакопеременной, линейной или унитарной группе.
 \end{Theo}

\section{Предварительные результаты}

\subsection{Редукция к почти простым группам и общие леммы}

Следуя \cite{BS_odd}, используем следующее обозначение.  Пусть $\pi$  --- некоторое множество простых
чисел, $m$~--- неотрицательное целое число.   Для конечной группы $G$ будем писать $G\in{{\mathcal B}{\mathcal S}}_{\pi}^{m}$, если группа $G$ обладает следующим свойством: класс сопряженности $D$ группы $G$ тогда и только тогда содержится в~$\Oo_\pi(G)$, когда любые  $m$ элементов из $D$ порождают $\pi$-группу.

Для подмножества $\pi$ множества $\mathbb{P}$ всех простых чисел полагаем $\pi'=\mathbb{P}\setminus\pi$.

\begin{Lemma}\label{red} {\rm \cite[лемма 7]{BS_odd}} {\it Пусть $\mathcal{X}$~--- класс конечных групп, замкнутый относительно взятия нормальных подгрупп, гомоморфных образов и расширений и содержащий все $\pi$-группы. Допустим,  $\mathcal{X}\nsubseteq\BS_{\pi}^{m}$ для некоторого натурального $m\geq 2$, и группа ${G\in\BS_{\pi}^{m}\setminus\mathcal{X}}$ выбрана так, что ее порядок является наименьшим. Тогда группа $G$ содержит подгруппу $L$ и элемент $x$ такие, что
\begin{itemize}
\item[$(1)$] $L\trianglelefteqslant G$;
\item[$(2)$] $L$ является неабелевой простой группой;
\item[$(3)$] $L$ не является $\pi$- или $\pi'$-группой;
\item[$(4)$] $C_G(L)=1$;
\item[$(5)$] любые $m$ сопряженных с $x$ элементов порождают $\pi$-группу;
\item[$(6)$] $x$ имеет простой порядок, принадлежащий $\pi$;
 \item[$(7)$] $G=\langle x, L\rangle$.
\end{itemize}}
\end{Lemma}

\begin{Lemma} \label{2notinpi} {\rm \cite[теорема 1]{BS_odd}} {\it Если $2\not\in\pi$, то $\BS_\pi^2$ совпадает с классом всех групп.
}
\end{Lemma}

\begin{Lemma} \label{guest} {\rm \cite[лемма~15]{GuestLevy}} {\it
 Пусть $G$~--- конечная группа и $x\in\Aut(G)$~--- автоморфизм, порядок которого равен степени простого числа $p$. Положим $M=C_G(x)$ и допустим, что $p$ делит $|G:M|$ и либо $M=N_G(M)$, либо $\mathrm{Z}(M) =1$. Тогда $x$ нормализует, но не централизует некоторую подгруппу в $G$, сопряженную с $M$.}
\end{Lemma}

Далее в индукционных рассуждениях будет использоваться следующие две очевидные леммы.

\begin{Lemma}\label{inductionLemma}
Пусть $L$~--- неабелева конечная простая группа порядка, кратного простому числу~$r$, и $x\in \Aut(L)^\sharp$. Предположим, что $x$ оставляет инвариантными некоторую подгруппу $H$ и её нормальную подгруппу $N$, индуцируя таким образом автоморфизм факторгруппы $\overline{H}=H/N$, который мы обозначим через $\overline{x}$. Предположим также, что $\overline{H}$ содержит простую $\overline{x}$-инвариантную подгруппу $\overline{L}$, на которой $\bar{x}$ действует нетождественно и порядок которой делится на некоторое простое число $r$.
Тогда $\beta_r(x,L)\leq \beta_r(\overline{\vphantom{P}x},\overline{L})$.
\end{Lemma}

\begin{proof}
Обозначим чертой $\overline{\phantom{G}}:\langle H,x\rangle\rightarrow \langle \overline{H},\overline{\vphantom{H}x}\rangle $ естественный гомоморфизм. Пусть $\beta_r(\bar{x},\overline{L})=k$. Это значит, что существуют $\bar{g}_1,\ldots,\bar{g}_k\in\overline{ L}$ такие, что порядок $\vert\langle \bar{x}^{\bar{g}_1},\ldots,\bar{x}^{\bar{g}_k}\rangle\vert$ делится на $r$. Пусть $g_1,\ldots,g_k$~--- какие-нибудь прообразы элементов $\bar{g}_1,\ldots,\bar{g}_k$. Тогда $$\overline{\langle x^{g_1},\ldots,x^{g_k}\rangle}=\langle \bar{x}^{\bar{g}_1},\ldots,\bar{x}^{\bar{g}_k}\rangle,$$ следовательно, $\vert \langle x^{g_1},\ldots,x^{g_k}\rangle\vert$ делится на $r$, откуда следует утверждение леммы.
\end{proof}

\begin{Lemma}\label{estim}
Пусть $L$~--- неабелева конечная простая группа, $x,y\in \Aut(L)^\sharp$. Предположим, что $x\in \langle y^{g_1},\dots,y^{g_k}\rangle$ для некоторых $g_1,\dots,g_k\in L$. Тогда
 $$
 \beta_r(x,L)\leq k\cdot\beta_r(y,L)
 $$
 для любого простого делителя $r$ порядка группы $L$.
\end{Lemma}

\subsection{Обозначения и предварительные сведения,\\ касающиеся классических групп}

 Напомним, что для векторного пространства $V$ символ $V^\sharp$ обозначает множество всех ненулевых векторов этого пространства. Везде далее предполагается, что $q=p^k$ для некоторого простого числа $p$ и $k\geq1$, через $\F_q$ обозначено конечное поле порядка~$q$.

Для вектора $v\in V$ и элемента $g\in \GL(V)$ образ $v$ под действием $g$ обозначен через $vg$. Если $x,g\in G$, то $x^g=g^{-1}xg$ и $[x,g]=x^{-1}x^g$. Если $v\in V$ и $g\in \GL(V)$, то $[v,g]=vg-v$. Для подгруппы $G$ группы $\GL(V)$ определим $[v,G]=\langle [v,g]\mid g\in G\rangle$. Заметим, что, поскольку $[v,G]=\langle vg-vh\mid g,h\in G\rangle$, $$\text{подпространство }[v,G] \text{ пространства }V \text{ является }G\text{-ин\-ва\-ри\-ант\-ным}.$$

Для группы $G$ через $\Oo_\infty(G)$, $\mathrm{Z}(G)$, $F(G)$, $\mathrm{F}^*(G)$, $\Phi(G)$, $G'=[G,G]$ и $G^{\infty}=\Oo^{\infty}(G)$ обозначены разрешимый радикал, центр, подгруппа Фиттинга, обобщённая подгруппа Фиттинга, подгруппа Фраттини коммутант и последний член ряда коммутантов (разрешимый корадикал) группы $G$ соответственно. Для любого множества $\pi$ простых чисел $\Oo^\pi(G)$~--- $\pi$-корадикал группы $G$, т.\,е. наименьшая по включению нормальная подгруппа, факторгруппа по которой является~$\pi$-группой. Как сказано выше, для множества $
\pi$ простых чисел его дополнение во множестве всех простых чисел обозначается через $\pi'$, поэтому группу $\Oo^\pi(G)$ можно трактовать также, как подгруппу, порожденную всеми $\pi'$-подгруппами.

Если пространство $V$ снабжено невырожденной билинейной, эрмитовой или квадратичной формой, то мы будем говорить, что подпространство $U$ пространства $V$ {\em невырождено}, если ограничение соответствующей формы на $U$ невырождено, и будем говорить, что  $U$ {\em вполне изотропно}, если ограничение на него соответствующей формы тождественно нулевое.

Мы всегда считаем, что $q$~--- степень простого числа~$p$. Образ в группе $\mathrm{PGL}_n(q)$ матрицы $\left(a_{ij}\right)\in\GL_n(q)$, заданной своими элементами, будем обозначать с помощью тех же элементов, заключенных в квадратные скобки~$\left[a_{ij}\right]$.

Далее через $\tau$ всегда будем обозначать автоморфизм группы $\GL_n(q)$, действующий по правилу $$\tau:A\mapsto (A^{-1})^\top,$$ где ${}^\top$~--- символ транспонирования матрицы. Через $\varphi_{p^m}$ обозначен автоморфизм группы $\GL_n(q)$, действующий по правилу $$\varphi_{p^m}:(a_{ij})\mapsto (a_{ij}^{p^m}).$$ Теми же символами $\tau$ и $\varphi_{p^m}$ мы обозначаем индуцированные $\tau$ и $\varphi_{p^m}$ автоморфизмы групп  $\SL_n(q)$, $\mathrm{PGL}_n(q)$ $\mathrm{PGL}_n(q)$ и $\mathrm{PSL}_n(q)$. В частности, если $q=p^k$ и $r$ делит $k$, то $\varphi_{q^{1/r}}$~--- полевой автоморфизм порядка~$r$.

 Мы как правило отождествляем группу $\mathrm{PSU}_n(q)$ с $\Oo^{p'}\left(C_{\mathrm{PGL}_n(q^2)}\left(\tau\varphi_{q^{1/2}}\right)\right)$. Как обычно, при изучении линейных и унитарных групп мы будем часто использовать унифицированное обозначение $\mathrm{PSL}_n^\varepsilon(q)$, где $\varepsilon=\pm1$,  а также знак соответствующего числа, полагая $$\mathrm{PSL}_n^{+}(q)=\mathrm{PSL}_n(q)
 \text{ и }\mathrm{PSL}_n^{-}(q)=\mathrm{PSU}_n(q).$$ Аналогичным образом определены $\mathrm{PGL}^\varepsilon_n(q)$, $\mathrm{GL}^\varepsilon_n(q)$ и $\mathrm{SL}^\varepsilon_n(q)$. Следуя \cite{atlas}, мы используем также краткую нотацию $$L_n(q)=L_n^+(q)=\mathrm{PSL}_n(q) \text{ и }U_n(q)=L_n^-(q)=\mathrm{PSU}_n(q).$$ При использовании лиевской нотации, мы считаем
 $$A_{n-1}(q)=A_{n-1}^+(q)=\mathrm{PSL}_n(q) \text{ и }{}^2A_{n-1}(q)=A_{n-1}^-(q)=\mathrm{PSU}_n(q).$$
 Символом $E_k$ обозначается единичная $(k\times k)$-матрица, а символом $A\otimes B$~--- кронекерово произведение матриц.

Применительно к автоморфизмам   групп лиева типа мы будем использовать следующую терминологию, близкую к принятой в \cite{GLS3} и несколько отличающуюся от принятой в \cite{Car}.
Понятия внутренне-диагонального автоморфизма в \cite{Car}
и \cite{GLS3} совпадают и не отличаются от используемого нами. В \cite[определение 2.5.10]{GLS3} введены подгруппы $\Phi_K$ и $\Gamma_K$ в группе автоморфизмов произвольной группы лиева типа $K$. Для групп лиева типа мы будем использовать обычно букву $L$, поэтому соответствующие подгруппы будем обозначать через $\Phi_L$ и $\Gamma_L$. Эти подгруппы можно отождествить, соответственно, с группами полевых и графовых автоморфизмов группы $L$ в смысле \cite{Car}. Через $\L$ обозначается группа внутренне-диагональных автоморфизмов группы~$L$.

  Применительно к линейным и унитарным группам скажем об автоморфизмах более подробно. В случае, когда $L=\PSL^\varepsilon_n(q)$ имеем $\L=\mathrm{PGL}^\varepsilon_n(q)$.

\begin{Lemma}\label{Aut} {\em \cite[теорема 2.5.12]{GLS3}} Пусть $L$ --- простая линейная или унитарная группа над полем ${\mathbb F}_q$ характеристики $p$. Тогда $\Aut(L)$ является расщепляемым расширением группы $\L$ с помощью абелевой группы $\Phi_L\Gamma_L$. При этом $\Phi_L\Gamma_L\cong\Phi_L\times\Gamma_L$.
\end{Lemma}

Для линейной простой группы $L$
автоморфизм $x\in\Aut(L)\setminus\L$ будем называть
\begin{itemize}
\item[]{\it полевым по модулю $\L$}, если 
образ $x$ в $\Aut(L)/\L$ лежит в $\L\Phi_L/\L$; при этом элементы группы $\Phi_L$ мы будем называть {\it каноническими полевыми автоморфизмами} группы $L$;
 \item[] {\it графовым по модулю} $\L$, если
образ $x$ в $\Aut(L)/\L$ лежит в $\L\Gamma_L/\L$; при этом элементы группы $\Gamma_L$ мы будем называть {\it каноническими графовым автоморфизмами} группы $L$;
\item[] {\it графово-полевым по модулю} $\L$ во всех остальных случаях; при этом  элементы  из $\Phi_L\Gamma_L\setminus(\Phi_L\cup\Gamma_L)$   будем называть {\it каноническими графово-полевыми автоморфизмами} группы $L$.
\end{itemize}

Пусть $L$ --- унитарная простая группа. Пусть $x\in\Aut(L)\setminus\L$. Автоморфизм $x$ будем называть
\begin{itemize}
\item[]{\it полевым по модулю $\L$}, если  порядок образа $x$ не делится на $2$; при этом такие элементы из группы $\Phi_L$ мы будем называть {\it каноническими полевыми автоморфизмами} группы $L$;
 \item[] {\it графовым по модулю} $\L$, если порядок образа $x$ равен  $2$; при этом такие элементы из группы $\Phi_L$ мы будем называть {\it каноническими графовыми автоморфизмами} группы $L$.
 \item[] {\it графово-полевым по модулю} $\L$, если порядок образа $x$ делится на $2$, но не равен $2$; при этом такие элементы из группы $\Phi_L$ мы будем называть {\it каноническими графово-полевыми автоморфизмами} группы $L$.

\end{itemize}

Заметим, что введенные понятия полевого и графово-полевого по модулю $\L$ автоморфизма $x$ группы $L$ совпадает с понятием полевого и графово-полевого автоморфизма в \cite[определение 2.5.13]{GLS3} соответственно в случае, когда $\langle x\rangle\cap \L=1$ (в частности, когда $x$ имеет простой порядок).

 \begin{center}\begin{table}
\caption
{Оценки на $\alpha_L(x)$ для классических групп $L$}\label{tab2}
\begin{center}\begin{tabular}{|c|c|c|c|r|}\hline
    \multirow{2}{*}{ $L$ }   & условия & условия & условия на  & \multirow{2}{*}{$\alpha(x,L)$} \\
  & на $n$   & на $q$  &  автоморфизм $x$     &   \\
\hline\hline
\multirow{18}{*}{$A_{n-1}(q)\cong L_n(q)$} & \multirow{10}{*}{$n=2$}  & \multirow{3}{*}{$q\ne 5,9$} & $|x|> 2$      & $2$   \\
                                               \cline{4-5}
                         &        &         & полевой,  $|x|=2$       & $\leq 4$  \\
                                              \cline{4-5}
                         &        &             &\multirow{2}{*}{ не полевой, $|x|=2$}     & \multirow{2}{*}{$3$}\\
                                  \cline{3-3}
                         &        & \multirow{4}{*}{$q= 9$}  &     & \\
                                             \cline{4-5}
                         &        &        & полевой,   $|x|=2$    & $5$     \\
                                             \cline{4-5}
                         &        &        & $|x|=3$       & $3$     \\
                                             \cline{4-5}
                         &        &       & $|x|>3$       & $2$    \\
                                  \cline{3-5}
                         &        & \multirow{3}{*}{$q=5$}    &  $|x|>2$      & $2$  \\
                                             \cline{4-5}
                         &        &          & не диагональный, $|x|=2$   & $3$  \\
                                             \cline{4-5}
                         &        &          & диагональный,  $|x|=2$      & $4$  \\
                                  \cline{3-5}

                         \cline{2-5}
                         & \multirow{3}{*}{ $n=3$} &          & не графово-полевой      & \multirow{2}{*}{$\leq 3$} \\
                         &        &          & или $|x|\ne 2$       &   \\
                                             \cline{4-5}
                         &        &          &    графово-полевой,  $|x|=2$     & $\leq 4$  \\
                         \cline{2-5}
                         &  \multirow{4}{*}{$n=4$} &  \multirow{2}{*}{$q>2$}  &графовый     & $\leq 6$  \\
                                            \cline{4-5}
                         &        &        &       \multirow{2}{*}{не графовый}      &  \multirow{2}{*}{$\leq 4$}  \\
                                  \cline{3-3}
                         &        & \multirow{2}{*}{$q=2$}  &      & \\
                                           \cline{4-5}
                         &        &        &    графовый     & $7$  \\
                         \cline{2-5}
                         &  \multirow{1}{*}{$n>4$} &          &               & $\leq n$  \\
                                  \hline   \hline
\multirow{11}{*}{${^2A}_{n-1}(q)\cong U_n(q)$} & \multirow{4}{*}{$n=3$}      &$q>3$    &                      &    $\leq 3$  \\
                                          \cline{3-5}
                             &            &\multirow{3}{*}{$q=3$}    & не внутренний        &    \multirow{2}{*}{$\leq 3$}  \\
                             &            &         & или $|x|\ne 2$       &              \\
                                                    \cline{4-5}
                             &            &         &внутренний, $|x|=2$   & $4$           \\
                             \cline{2-5}
                             & \multirow{6}{*}{$n=4$}      &  \multirow{2}{*}{$q>2$}   & не графовый         & $\leq 4$  \\
                                                     \cline{4-5}
                             &           &           &     \multirow{2}{*}{графовый }        &  \multirow{2}{*}{$\leq 6$ } \\
                                         \cline{3-3}
                             &           & \multirow{4}{*}{$q=2$}      &           &   \\  \cline{4-5}
                             &           &        & трансвекция         & $\leq 5$       \\
                                                    \cline{4-5}
                             &           &           & не трансвекция      & \multirow{2}{*}{$\leq 4$}  \\

                             &           &       & или не графовый     &      \\
                             \cline{2-5}
                             &  $n>4$    &           &                    & $\leq n$  \\
                                  \hline   \hline
\multirow{7}{*}{$C_{n}(q)\cong S_{2n}(q)$} &  \multirow{4}{*}{ $n=2$} &    \multirow{2}{*}{ $q>3$}   & $|x|= 2$           & $\leq 5$  \\
                                                     \cline{4-5}
                             &           &             &  \multirow{2}{*}{$|x|> 2$ }          &  \multirow{2}{*}{$\leq 4$ }   \\
                                         \cline{3-3}
                             &           & \multirow{2}{*}{$q=3$}     &           &  \\
                                                     \cline{4-5}
                             &           &           & $|x|= 2$           & $\leq 6$   \\
                             \cline{2-5}
                             &  \multirow{3}{*}{$n>2$}    &           &  не трансвекция    & $\leq n+3$  \\
                                         \cline{3-5}
                             &           & нечетно   &    \multirow{2}{*}{  трансвекция  }  & $\leq 2n$  \\
                                         \cline{3-3} \cline{5-5}
                             &           & четно     &       & $\leq 2n+1$  \\
                           \hline   \hline
 \multirow{4}{*}{$B_{n}(q)\cong O_{2n+1}(q)$} & \multirow{4}{*}{$n\geq 3$}  & \multirow{2}{*}{нечетно} &   отражение      &  $ 2n+1$    \\
                                                    \cline{4-5}
                             &            &         &  не  отражение        & \multirow{2}{*}{$\leq n+3$} \\
                                          \cline{3-4}
                             &            & \multirow{2}{*}{четно}   &    не трансвекция    &  \\
                                                     \cline{4-5}
                             &            &         &     трансвекция    & $2n+1$  \\
                             \hline   \hline

  \multirow{2}{*}{$D_{n}(q)\cong O_{2n}^+(q)$,} & \multirow{4}{*}{$n\geq 4$}  & \multirow{2}{*}{нечетно} &   отражение      &  $ 2n$    \\
                                                    \cline{4-5}
                             &            &         &  не  отражение        & \multirow{2}{*}{$\leq n+3$} \\
                                          \cline{3-4}
  \multirow{2}{*}{${}^2D_{n}(q)\cong O_{2n}^-(q)$}  &            & \multirow{2}{*}{четно}   &    не трансвекция    &  \\
                                                     \cline{4-5}
                             &            &         &     трансвекция    & $2n$  \\
                             \hline
\end{tabular}\end{center}
\end{table}\end{center}

Следующая лемма содержит информацию о параметре $\alpha(x,L)$, где $x$~--- автоморфизм простого порядка классической простой группы~$L$, взятую из~\cite{GS}.

  \begin{Lemma} \label{alpha_classic}
Пусть  $S$~--- простая классическая группа  и $x\in \Aut L$~---  элемент простого порядка. Тогда $\alpha(x,L)$ удовлетворяет условию, указанному в последнем столбце
таблицы~{\rm\ref{tab2}}.
\end{Lemma}

\begin{Lemma} \label{alpha_irreducible}
Пусть  $L$~--- простая классическая группа  и $x\in \widehat{L}$~--- нетривиальный элемент простого порядка индуцированный неприводимым полупростым
элементом соответствующей $L$ группы подобий. Тогда $\alpha(x,L)\leq 3$.
\end{Lemma}
\begin{proof}
 См. \cite[общее доказательство теорем 4.1--4.4; более точно, для линейных групп~--- стр.~534, для унитарных --- стр.~536, для симплектических~---
стр.~538 и для ортогональных~--- стр.~539]{GS}.
\end{proof}

\begin{Lemma}\label{betaL2q} Пусть $L=L_2(q)$ и $r\in \pi(L)$~--- нечетное число, не делящее $q$. Тогда для любой инволюции $x\in\L$ выполнено $\beta_r(x,L)=2$.
\end{Lemma}

\begin{proof} Из условия следует, что $r$ делит $q-\varepsilon$ для некоторого $\varepsilon\in\{+1,-1\}$. Допустим сначала, что $x\in L$. Поскольку все инволюции в $L$ сопряжены (если $q$ четно, это следует из теоремы Силова, равенства $L=\L$ и хорошо известного факта, что нормализатор силовской 2-подгруппы в $L$ является группой Фробениуса с ядром порядка $q$ и дополнением порядка $q-1$, регулярно действующим на неединичных 2-элементах силовской 2-подгруппы,  а если $q$ нечетно, см.~\cite[теорема 4.5.1 и таблица 4.5.1]{GLS3}), мы можем считать, что $x$ содержится в подгруппе диэдра $D$ порядка $q-\varepsilon$ и инвертирует все элементы единственной подгруппы $\langle y\rangle$ порядка $r$ в группе $D$. Тогда элемент $$xx^y=(y^{-1})^xy=y^2$$ имеет порядок~$r$ и $\beta_r(x,L)=2$, как и утверждается.

Если $x\notin L$, то $q$ нечетно. Все инволюции в $\L\setminus L$ сопряжены элементами из $L$, как следует из  \cite[теорема 4.5.3 и таблица 4.5.2]{GLS3}, причем сопряжены с инволюцией из подгруппы диэдра $\widehat{D}$ порядка $2(q-\varepsilon)$, которая инвертирует элементы единственной подгруппы  $\langle y\rangle$ порядка $r$ в группе $\widehat{D}$. Рассуждая как и в предыдущем случае, получаем требуемое.
\end{proof}

\begin{Lemma}\label{generationSL2p} Пусть $\F_q$~--- поле нечетной характеристики $p$ и $\beta\in\F^*_q$. Обозначим через $\F_{q^{\phantom{0}}_0}$ подполе в $\F_q$, порожденное $\beta^2$. Рассмотрим в $\SL_2(q)$ элементы
$$
x=\left(
\begin{array}{rr}
  1 &  \\
 \beta & 1
\end{array}
\right),
\quad
y=\left(
\begin{array}{rr}
  1 & \beta \\
 & 1
\end{array}
\right).
$$
 Тогда подгруппа
$
H=
\left\langle x,y
  \right\rangle$ группы $\SL_2(q)$ либо
 изоморфна $\SL_2(q_0)$, либо $q_0=9$ и $H$ изоморфна подгруппе в $\SL_2(q_0)$, образ которой в $\PSL_2(q_0)$ изоморфен $A_5$, а сама $H$ содержит подгруппу, изоморфную $\SL_2(3)$.  В частности, $H$ всегда содержит подгруппу, изоморфную $\SL_2(p)$, и, в частности, содержит содержит матрицу $$\left(
\begin{array}{rr}
  -1 &  \\
 & -1
\end{array}
\right)
,$$
которая является единственным элементом порядка $2$  группы $\SL_2(q)$.

\end{Lemma}

\begin{proof}
Рассмотрим матрицу
$$
g= \left(
\begin{array}{rr}
  1 &  \\
 & \beta
\end{array}
\right)\in \GL_2(q)
$$
и подгруппу $H^g=\langle x^g,y^g\rangle\cong H$. Непосредственными вычислениями убеждаемся, что
$$
x^g=g^{-1}xg=\left(
\begin{array}{cc}
  1 &  \\
 1 & 1
\end{array}
\right),
\quad
y^g=g^{-1}yg=\left(
\begin{array}{cc}
  1 & \beta^2 \\
 & 1
\end{array}
\right).
$$
Теперь из  \cite[гл.~2, теорема~8.4]{Gor} следует, что подгруппа $H^g$, а значит и $H$, такая, как утверждается в лемме.
\end{proof}

\begin{Lemma} \label{cases}
Пусть $\Delta=\Delta(V)$~--- группа всех подобий конечномерного векторного пространства $V$, снабженного невырожденной или тривиальной билинейной или эрмитовой формой. Пусть $x\in\Delta$~--- примарный элемент, образ которого в группе $\Delta(V)/\mathrm{Z}(\Delta(V))$ имеет простой порядок. Тогда имеет место один из следующих случаев:
\begin{itemize}
  \item[$(1)$] элемент $x$ унипотентен и стабилизирует подпространство размерности~$1$;
  \item[$(2)$] форма тривиальна, элемент $x$ полупрост и оставляет инвариантными два дополняющих друг друга ненулевых подпространства;
  \item[$(3)$] форма невырождена, элемент $x$ полупрост и оставляет инвариантными собственное ненулевое невырожденное подпространство и ортогональное дополнение к нему;
  \item[$(4)$] форма невырождена, имеет максимальный индекс Витта, элемент $x$ полупрост и оставляет инвариантным вполне изотропное подпространство максимальной размерности $\displaystyle\frac{1}{2}\dim V$;
  \item[$(5)$] $x$ полупрост и действует неприводимо на $V$.
\end{itemize}
\end{Lemma}

\begin{proof}
 Так как образ элемента $x$ в $\Delta(V)/\mathrm{Z}(\Delta(V))$~--- элемент простого порядка, $x$ либо унипотентен, либо полупрост. В первом случае пусть $k$~--- наименьшее натуральное число такое, что $(x-1)^k=0$. Тогда $(x-1)^{k-1}\ne 0$ и существует  вектор $v\in V$ такой, что $v(x-1)^{k-1}\ne 0$, а $v(x-1)^{k}= 0$. Таким образом,
для ненулевого вектора $u=v(x-1)^{k-1}$ выполнено $ux=u$ и $\langle u\rangle$~--- одномерное $x$-инвариантное подпространство, т.е. верно~$(1)$. 

Пусть теперь элемент $x$ полупрост. Если он неприводим, то верно~$(5)$. Если форма тривиальна, то по теореме Машке \cite[теорема~(1.9)]{Isaacs} верно $(2)$. Поэтому считаем, что $x$ приводим и форма невырождена. Пусть $U$~--- собственное ненулевое $x$-инвариантное подпространство. Элемент то $x$ стабилизирует ортогональное дополнение $U^\perp$. Если $U$ невырождено, то $U\cap U^\perp=0$ и $U^\perp$~--- дополнение в $U$ в~$V$, тем самым верно~$(3)$.

Если $U$ не является невырожденным, то его радикал $\mathrm{rad}\,U=U\cap U^\perp$~--- собственное ненулевое вполне изотропное $x$-инвариантное подпространство. Можно, таким образом, считать $U$ вполне изотропным. Если $\dim U<\displaystyle\frac{1}{2}\dim V$, то $U<U^\perp$, $U=\mathrm{rad}\,U^\perp$ и по теореме Машке $U^\perp$ обладает $x$-инвариантным дополнением $W$ к~$U$, причем $W\cong U^\perp/U$ невырождено, т.е. имеет место~$(3)$. Если же $\dim U=\displaystyle\frac{1}{2}\dim V$, то выполнено~$(4)$.
\end{proof}

\begin{Lemma} \label{Parabolic} {\rm \cite[лемма 2.2]{GS}} { \it
Пусть
 $L$ --- простая группа лиева типа,  $G = \L$ и $x\in G$. Тогда выполнены следующие утверждения.
\begin{itemize}
\item[$(1)$] Если элемент  $x$ унипотентный, то пусть $P_1$ и $P_2$~--- различные параболические максимальные подгруппы в~$G$, содержащие общую подгруппу Бореля, и  $U_1$,  $U_2$~--- унипотентные радикалы подгрупп $P_1$ и $P_2$ соответственно. Тогда элемент $x$ сопряжен с элементом из $P_i\setminus U_i$ для $i = 1$ или $i = 2$.
\item[$(2)$]  Если элемент $x$ полупростой, допустим, что $x$ содержится в некоторой параболической подгруппе группы~$G$. Если ранг группы $L$ не меньше двух, то существует параболическая максимальная подгруппа~$P$ c дополнением Леви $J$ такая, что элемент $x$ сопряжен с некоторым элементом из $J$, не централизующим никакую из (возможно, разрешимых) компонент Леви группы~$J$.
\end{itemize}}
\end{Lemma}

\begin{Lemma} \label{UnipotentUnitary} Пусть $V$~--- пространство с невырожденной эрмитовой формой над полем $\F_{q^2}$ нечетной размерности $\dim V=n\geq 5$. Тогда для унипотентного элемента $x\in\SU(W)$ простого порядка имеет место один из следующих случаев:
\begin{itemize}
  \item[$(1)$] $x$  стабилизирует невырожденное подпространство размерности~$1$;
  \item[$(2)$] $x$ стабилизирует максимальное вполне изотропное подпространство и индуцирует на нем нетождественное преобразование;
  \item[$(3)$] $q$ нечетно и существует сопряженный с $x$ элемент $x^g$ такой, что в подгруппе $\langle x, x^g\rangle$ некоторая инволюция стабилизирует максимальное вполне изотропное подпространство и индуцирует на нем нескалярное преобразование.
\end{itemize}
\end{Lemma}

\begin{proof} Для $\alpha\in\F_{q^2}$ положим $\bar{\alpha}=\alpha^q$.

 В~\cite[предложение~2.2]{GLO'B} найдены представители классов сопряженности унипотентных элементов группы ${\SU(V)\cong\SU_n(q)}$. Известно, что с точностью до сопряженности элементом из $\GU(V)$ унипотентному элементу $x$ соответствует однозначно определенное неупорядоченное разбиение
$$n=n_1+\dots+n_s$$ и соответствующее ему разложение
$$
V=V_1\oplus\dots\oplus V_s
$$
 пространства $V$ в ортогональную сумму невырожденных $x$-инвариантных подпространств $V_1,\dots,V_s$ таких, что $\dim V_t=n_t$, для $t\in\{1,\dots,s\}$, и действие $x$ на $V_t$ задается следующим образом.

 Если $n_t=2k$, $k=k(t)$, то в $V_t$ существуют упорядоченный  базис
 $$e^t_1,\dots,e^t_k,f^t_k,\dots,f^t_1$$
 такой, что
 $$(e_i^t,e_j^t)=(f_i^t,f_j^t)=0, \quad(e_i^t,f_j^t)=\delta_{ij}\quad\text{для всех }i,j=1\dots,k,$$
 и элемент $\beta\in\F_{q^2}^*$ такой, что $\beta+\bar{\beta}=0$, для которых выполнены равенства
 \begin{eqnarray*}
   e^t_ix= & e^t_i+\dots+e^t_k+\beta f^t_k &\text{ для всех } i=1,\dots,k, \\
   f^t_ix=
   & f^t_i- f^t_{i-1} & \text{ для всех } i=2,\dots,k, \\
   f^t_1x= & \phantom{f^t_i- }   f^t_{1}.&
  \end{eqnarray*}

  Если же $n_t=2k+1$, $k=k(t)$, то в $V_t$ существуют упорядоченный  базис
 $$e^t_1,\dots,e^t_k,d^t,f^t_k,\dots,f^t_1,$$
 такой, что
 $$(e_i^t,e_j^t)=(f_i^t,f_j^t)=(e_i^t,d^t)=(f_i^t,d^t)=0,$$
 $$(d^t,d^t)=1, (e_i^t,f_j^t)=\delta_{ij}\text{ для всех }i,j=1\dots,k,$$
 и элемент $\gamma\in\F_{q^2}^*$ такого, что $\gamma+\bar{\gamma}=-1$, для которых выполнены равенства
 \begin{eqnarray*}
   e^t_ix= & e^t_i+\dots+e^t_k+d^t+\gamma f^t_k &\text{ для всех } i=1,\dots,k, \\
   d^tx=
   & d^t- f^t_{k\phantom{-1}} &\\
   f^t_ix=
   & f^t_i- f^t_{i-1} & \text{ для всех } i=2,\dots,k, \\
   f^t_1x= & \phantom{f^t_i- }   f^t_{1}.&
  \end{eqnarray*}
Ясно, что случай $(1)$, т.\,е. существование одномерного $x$-инвариантного невырожденного подпространства, равносилен тому, что  $n_t=1$ для некоторого $t\in\{1,\dots,s\}$. Поэтому считаем, что $n_t>1$ для всех~$t$. Поскольку $n$ нечетно, некоторое $n_t$ также нечетно и для него $n_t\geq 3$. Теперь из того, что $x$ имеет порядок $p$, заключаем, что $p>2$. В самом деле, при $p=2$ элемент $x^2$ нетождественно действует на $V_t$:
$$e^t_kx^2=(e^t_k+d^t+\gamma f_k^t)x=e_k^t+d^t+\gamma f_k^t+d^t+f_k^t+\gamma(f_k^t+f)=e_k^t+f_k^t+\gamma f\ne e_k^t,$$  где  $f=0\text{ или } f=f^t_{k-1}.$

Обозначим через $U$ подпространство, порожденное всеми $f^t_i$, где $t=1,\dots,s$, $i=1,\dots, k(t)$. Ясно, что $U$ вполне изотропно и $x$-инвариантно.

 Предположим, что  ровно одно слагаемое  в сумме $n_1+\dots+n_s$ нечетно.Тогда подпространство $U$ является максимальным вполне изотропным. При этом если $k=k(t)>1$ для некоторого $t\in\{1,\dots,s\}$, то из равенства
$$
f^t_{k}x=f^t_k-f^t_{k-1}
$$
следует, что $x$ действует на $U$ нетождественно, т.\,е. имеет место случай $(2)$ леммы. Поэтому мы можем считать, что $$n_1,\dots,n_{s-1}=2 \quad\text{и}\quad n_s=3.$$ Так как $n\geq 5$, видим, что $s\geq 2$ и $n_1=2$. Рассмотрим элемент $g\in \SL(V)$ такой, что
$$e^1_1g=f^1_1,\quad f^1_1g=-e_1^1,$$
$$e^t_1g=e^t_1,\quad f^t_1g=f^t_1\quad\text{при } t>1\text{ и}$$
$$d^sg=d^s.$$  Тогда матрицы элементов $x$ и $x^g$ в базисе
$$
e^1_1,f^1_1,e^2_1,f^2_1,\dots, e^{s-1}_1,f^{s-1}_1,e^{s}_1,d^s,f^{s}_1
$$
имеют вид
$$
\left(
\begin{array}{rrrrrrrrrr}
  1 & \beta &  &  &  &  &  &  &  &  \\
   & 1 &  &  &  &  &  &  &  &  \\
   &  & 1 & \beta &  &  &  &  &  &  \\
   &  &  & 1 &  &  &  &  &  &  \\
   &  &  &  & \ddots &  &  &  &  &  \\
   &  &  &  &  & 1 & \beta &  &  &  \\
   &  &  &  &  &  & 1 &  &  &  \\
   &  &  &  &  &  &  & 1 & 1 & \gamma \\
   &  &  &  &  &  &  &  & 1 & -1 \\
   &  &  &  &  &  &  &  &  & 1
\end{array}
\right)\quad\text{и}\quad
\left(
\begin{array}{rrrrrrrrrr}
  1 & &  &  &  &  &  &  &  &  \\
  -\beta & 1 &  &  &  &  &  &  &  &  \\
   &  & 1 & \beta &  &  &  &  &  &  \\
   &  &  & 1 &  &  &  &  &  &  \\
   &  &  &  & \ddots &  &  &  &  &  \\
   &  &  &  &  & 1 & \beta &  &  &  \\
   &  &  &  &  &  & 1 &  &  &  \\
   &  &  &  &  &  &  & 1 & 1 & \gamma \\
   &  &  &  &  &  &  &  & 1 & -1 \\
   &  &  &  &  &  &  &  &  & 1
\end{array}
\right)
$$
 соответственно, где элементы $\beta$ и $\gamma$ из $\F_{q^2}$ выбраны описанным выше способом. Далее, рассмотрим подгруппу $$H=\langle x,x^g\rangle=\left\langle \left(x^g\right)^{-1},x\right\rangle.$$ Ее можно отождествить с некоторой подгруппой в группе $K$ всех блочно-диа\-го\-наль\-ных матриц вида
  $$
  \left(
\begin{array}{rr}
  A &  \\
  & B
\end{array}
\right),
  $$
  где $A\in\SL_2(q^2)$, a $B$~--- произвольная унитреугольная матрица над $\F_{q^2}$ порядка~${n-2}$. Рассмотрим также   эпиморфизм
$$
\overline{\phantom{x}\vphantom{\left(x^g\right)}}:K\rightarrow \SL_2(q^2),
$$
действующий по правилу
$$
  \left(
\begin{array}{rr}
  A &  \\
  & B
\end{array}
\right)\mapsto A.
  $$
Тогда
$$
\overline{\left(x^g\right)}^{-1}=\left(
\begin{array}{rr}
  1 &  \\
 \beta & 1
\end{array}
\right)\quad\text{и}\quad
\overline{x\vphantom{\left(x^g\right)}}=\left(
\begin{array}{rr}
  1 & \beta \\
 & 1
\end{array}
\right).
$$
Поэтому, как следует из леммы~\ref{generationSL2p}, подгруппа
$$
\overline{H\vphantom{\left(x^g\right)^{-1}}}= \overline{\left\langle
\left(x^g\right)^{-1},{x}
  \right\rangle}\leq\overline{K\vphantom{\left(x^g\right)^{-1}}}=\SL_2(q^{{2}}).
$$
содержит единственный элемент $\overline{y}$ порядка~2 группы $SL_2(q^{{2}})$, и это   матрица
$$
\left(
\begin{array}{rr}
  -1 &  \\
 & -1
\end{array}
\right).$$
Возьмем  его прообраз $y$ в $H$ также порядка~$2$. Тогда матрица $y$ в базисе
$$
e^1_1,f^1_1,e^2_1,f^2_1,\dots, e^{s-1}_1,f^{s-1}_1,e^{s}_1,d^s,f^{s}_1
$$  имеет вид\footnote{Можно заметить, что на месте символов $*$ стоят нули, так как в противном случае порядок образа элемента $y$ в некотором гомоморфном образе подгруппы $H$ был бы кратен $p$, вопреки выбору~$y$, но для наших рассуждений это не важно.}
$$
\left(
\begin{array}{rrrrrrrrrr}
  -1 & &  &  &  &  &  &  &  &  \\
   & -1 &  &  &  &  &  &  &  &  \\
   &  & 1 & * &  &  &  &  &  &  \\
   &  &  & 1 &  &  &  &  &  &  \\
   &  &  &  & \ddots &  &  &  &  &  \\
   &  &  &  &  & 1 & * &  &  &  \\
   &  &  &  &  &  & 1 &  &  &  \\
   &  &  &  &  &  &  & 1 & * & * \\
   &  &  &  &  &  &  &  & 1 & * \\
   &  &  &  &  &  &  &  &  & 1
\end{array}
\right),
$$
и векторы $f_1^1,f_1^2,\dots f_1^{s-1}, f^s_1$ будут собственными для $y$: $$f_1^1y=-f^1_1\quad\text{и}\quad f_1^ty=f_1^t\quad\text{при}\quad t>1.$$
 Следовательно, максимальное вполне изотропное подпространство $$U=\langle f_1^1,\dots f_1^{s-1}, f^s_1\rangle$$ инвариантно относительно элемента $y$, и на $U$ этот элемент индуцирует нескалярное преобразование. Таким образом, имеет место случай~$(3)$.

Покажем, наконец, что если среди $n_1,\dots, n_s$ более одного нечетного числа, то снова имеет место случай $(2)$. Не уменьшая общности, можем считать, что нечетными будут числа $n_1,\dots, n_{s'}$, а числа $n_{s'+1},\dots, n_s$ четны, причем, поскольку $n$ нечетно, число $s'$ также нечетно. Рассмотрим ортогональное дополнение $U^\perp$ к $U$. Оно $x$-инвариантно, и из того, что $\mathrm{codim}\, U^\perp=\dim U=n-\dim U^\perp$ совпадает с числом векторов $e_i^t$, $t=1,\dots,s$, $i=1,\dots, k_t$, видно, что $$U^\perp=\langle d^{t'}, f_i^t\mid t'=1,\dots,s',\quad t=1,\dots,s,\quad i=1,\dots, k_t\rangle.$$

Рассмотрим $(s'-1)/2$ пар векторов $(d^1,d^2),\dots,(d^{s'-2},d^{s'-1})$. Возьмем одну такую пару $(d^{t'},d^{t'+1})=(d,d')$ и построим по ней пару $(a,b)=(a^{t'},b^{t'})$ по следующему правилу. Возьмем $\mu\in \F_{q^2}$ так, чтобы $\mu\bar{\mu}=-1$, и положим
$$a=d+\mu d', \quad b=d-\mu d'.$$
Так как $(d,d')= 0$ и $(d,d)=(d',d')=1$, имеем  $$(a,a)=(d+\mu d',d+\mu d')=1+\mu\bar{\mu}=0$$ и аналогично $(b,b)=0$. Заметим, что $$dx+U= d+U\text{ и }d'x+U=d'+U,$$ т.е. $x$ индуцирует тождественное преобразование факторпространства  $$(\langle d,d'\rangle+U)/U=(\langle a,b\rangle+U)/U.$$
Теперь видим, что
$$W=\langle a^1,a^3,\dots,a^{s'-2}, U\rangle$$
является максимальным вполне изотропным подпространством пространства $V$, которое к тому же $x$-инвариантно. При этом $$a^1x=(d^1+\mu d^2)x=d^1+\mu d^2 - f^1_{k(1)}-\mu f^1_{k(2)}\ne d^1+\mu d^2= a^1,$$ т.\,е. $x$ действует на $W$ нетождественно, и имеет место случай~$(2)$.
\end{proof}

\begin{Lemma}\label{irredlinear}
 Пусть  имеет место одно из утверждений:
  \begin{itemize}
    \item[$(+)$] $V$~--- конечномерное векторное пространство над конечным полем $F$, $V=U\oplus W$ для некоторых подпространств $U$ и $W$ и $2\leq \dim U\leq \dim W$;
    \item[$(-)$] выполнено то же, что и в утверждении $(+)$ и, кроме того, на $V$ задана невырожденная эрмитова форма ограничения которой на $U$ и $W$ невырождены, причем $U\perp W$.
  \end{itemize}
  Пусть знак $\varepsilon\in\{+,-\}$ показывает, какое из данных утверждений имеет место, и пусть $F=\mathbb{F}_q$ при ${\varepsilon=+}$ и $F=\mathbb{F}_{q^2}$ при ${\varepsilon=-}$, а $n=\dim V$.

    Допустим, подгруппы $L^{\phantom{x}}_U\leq \GL^\varepsilon(U)$  и $L_W^{\phantom{x}}\leq \GL^\varepsilon(W)$ неприводимы,  и определим естественным образом $L=L^{\phantom{x}}_U\times L^{\phantom{x}}_W$ как подгруппу в~$\GL^\varepsilon(V)$.
Предположим, что  элемент $x\in L$ таков, что его естественная проекция на $L^{\phantom{x}}_U$ действует неприводимо на~$U$.

Тогда справедливы следующие утверждения.
\begin{itemize}
  \item[$(1)$] Cуществует элемент $g\in \SL^\varepsilon(V)$ такой, что подгруппа $G=\langle L,x^g\rangle$ неприводимо действует на $V$.
  \item[$(2)$] Предположим, что $\dim U\leq\dim W$ и $\SL^\varepsilon(W)\leq L^{\phantom{x}}_W\leq\GL^\varepsilon(W)$. Тогда либо подгруппа $G\leq
  \GL^\varepsilon(V)$, определенная в $(1)$, примитивна,  либо мощность системы импримитивности равна $n$ и $(n,q)\in\{(4,2),(4,3)\}$.
  \item[$(3)$] Допустим,  $(n,q)\notin\{(4,2),(4,3),(4,5)\}$, $|x|$~--- простое число и элемент $x$ при $\varepsilon=+$ не имеет собственных векторов, a при $\varepsilon=-$ не имеет невырожденных собственных векторов.  Пусть $U$ является $x$-инвариантным подпространством наименьшей размерности, $t=\mathrm{codim}\, U$. Положим $$m=\left\{\begin{array}{cc}
                                                                                                                          t & \text{ при } t>2, \\
                                                                                                                          3 &   \text{ при } t=2.
                                                                                                                        \end{array}
  \right.$$
  Тогда некоторые $m+1$ элементов, сопряженных с $x$ посредством элементов из $\SL^\varepsilon(V)$,  порождают подгруппу, содержащую одну из подгрупп
  \begin{itemize}
    \item[$\bullet$] $\SL^\varepsilon_n(q)$,
    \item[$\bullet$] $\Sp_n(q)$,
    \item[$\bullet$] $\SL_n^-(q_0)$ при $q=q_0^2$,
    \item[$\bullet$] $\Oo^\pm_n(q)$ при четных $q$,
    \item[$\bullet$] $\Sym_{n+1}$ при $q=2$ и нечетном $n> 6$ и
    \item[$\bullet$] $\Sym_{n+2}$ при $q=2$ и четном $n\geq 6$.
  \end{itemize}
\end{itemize}
\end{Lemma}

\begin{proof}
Покажем сначала, что
\begin{itemize}
  \item[$(*)$] $U$ и $W$~--- это единственные собственные $L$-инвариантные подпространства пространства $V$.
\end{itemize}
Действительно, предположим, что собственное подпространство $X$ пространства $V$ является $L$-инвариантным. Поскольку $2\leq \dim U\leq \dim W$ и подгруппы $L_U$ и $L_W$ неприводимы на соответствующих подпространствах, для любых $u\in U^\sharp $  и $w\in W^\sharp$ справедливы равенства $[u,L_U]=U$ и $[w,L_W]=W$. Пусть $v\in X^\sharp$. Тогда вектор $v$ единственным образом представим в виде $v=u+w$, где $u\in U$ и $w\in W$. Если $u\not=0$, то в силу неприводимости группы~$L_U$ имеем $$X\geq[ v,G] \geq \left[ v,L_U^{\phantom{x}}\right]=\left[u,L_U^{\phantom{x}}\right] =U.$$ Следовательно, если $u\not=0$, то $X$ содержит $U$. Аналогичные рассуждения показывают, что если $w\not=0$, то $X$ содержит $W$. Утверждение $(*)$ доказано.

Далее, мы утверждаем, что

\begin{itemize}
  \item[$(**)$] существует такой элемент $g\in \SL^\varepsilon(V)$, что
$$Ug\cap U\not=0\quad\text{ и }\quad Ug\cap W\not=0,$$ в частности $Ug\nleqslant U$ и $Ug\nleqslant W$.
\end{itemize}
Действительно, выберем в $U$ и $W$ базисы $e_1,\dots,e_s$ и $f_1,\dots,f_t$ соответственно,  которые в случае $(-)$ будут ортонормированными. Поскольку $s=\dim U\geq 2$, утверждение $(**)$ выполнено для элемента $g\in\SL^\varepsilon(V)$, определенного равенствами
$$ e_ig=e_i\text{ при } i<s, \quad  f_ig= f_i\text{ при } i<t,\quad  e_sg=- f_t,\quad \text{ и }\quad  f_tg= e_s.$$

Покажем, что подгруппа $\langle L,x^g\rangle$ неприводима. Допустим, $X$~--- собственное ненулевое $\langle L,x^g\rangle$-инвариантное подпространство. Тогда оно является $G$-ин\-ва\-ри\-ант\-ным, а значит либо $X=U$, либо $X=W$. С другой стороны $Ug$~--- $x^g$-не\-при\-во\-ди\-мое $x^g$-инвариантное подпространство, значит $Ug\cap X$ также является $x^g$-инвариантным, что противоречит $x^g$-неприводимости подпространства~$Ug$. Утверждение $(1)$ доказано.

Прежде, чем доказывать утверждение $(2)$, отметим, что пара подпространств $\{U,W\}$ не является системой импримитивности для группы~$G$. В самом деле, в противном случае из того, что $e_1x,e_1g\in U$ следовало бы, что $e_1x^g\in U$ и $x^g$ стабилизирует блоки этой системы, как и подгруппа $L$. Значит группа $G=\langle L,x^g\rangle$ была бы приводима, вопреки~$(1)$.

Предположим, что что
$$
V=V_1\oplus\dots\oplus V_m,
$$
$m>1$ и $V_i^y\in\Omega:=\{V_1,\dots, V_m\}$ для любых $i\in\{1,\dots,m\}$ и $y\in G$. Из неприводимости группы $G$ заключаем также, что $G$ действует транзитивно на $\Omega$, в частности все $V_i$ изометричны и их размерности одинаковы. Пусть $\dim V_i=k$.

Поскольку сумма элементов любой $L$-орбиты на  $\Omega$ инвариантна относительно $L$, из утверждения~$(1)$ заключаем, что либо $L$ действует транзитивно на $\Omega$, имеется ровно две $L$-орбиты на этом множестве, сумма элементов одной из которых равна $U$, а другой $W$.

Легко видеть, что случай нетранзитивного действия $L$ на $\Omega$ эквивалентен существованию $V_i\in\Omega$ такого, что $V_i\leq W$. Заметим, что в этом случае $\{V_i\mid V_i\leq W\}$~--- система импримитивности для групп $L$ и $L_W$ на~$W$.
Пусть $$\SL^\varepsilon(W)\leq L_W\leq\GL^\varepsilon(W),$$ как в утверждении~$(2)$.
Тогда группа $L_W$ примитивна, и значит $W$ совпадает с~$V_i$.
 При этом
  $$ U=\sum\limits_{j\ne i}V_j,$$ как следует из~$(*)$.
  Учитывая, что $$(m-1)k=\dim U\leq\dim W=\dim V_i=k,$$ получаем $m=2$ и $\Omega =\{U,W\}$. Но, как мы отметили выше, $\{U,W\}$ не является системой импримитивности для группы~$G$.

Рассмотрим случай, когда $L$ транзитивна на $\Omega$. Пусть
$$
\omega:V\rightarrow W
$$
обозначает проекцию на $W$ параллельно~$U$. Тогда $vy\omega=v\omega y$ для любых ${v\in V}$ и ${y\in L}$. В частности, $L_W$ транзитивно действует на $\Omega\omega=\{V_i\omega\mid i=1,\dots,m\}$ и размерности всех проекций $V_i\omega$ одинаковы.  Обозначим их через $k'$.

 Рассмотрим сначала случай~$(+)$.
  Группа $\SL(W)$ действует транзитивно на подпространствах пространства $W$ одинаковой размерности. Так как $\SL(W)\leq L_W$, множество $\Omega\omega$ совпадает с множеством $k'$-мерных подпространств пространства $W$, которых заведомо не меньше, чем одномерных подпространств в~$W$. Таким образом, если $\dim W=t$, то
$$ m=|\Omega|\geq |\Omega\omega|\geq \frac{q^t-1}{q-1}\geq 2^t-1\geq 2^{n/2}-1,$$ где $q$~--- порядок поля, над которым определено пространство~$V$.
Если $m<n$, то $m=n/k\leq  n/2$. Из условия следует, что $n\geq 4$, но при таких таких $n$ неравенство $n/2\geq 2^{n/2}-1$ неверно. Поэтому $m=n$. При $n>4$ нарушается неравенство
 $$
 n\geq 2^{t}-1, \text{ где } t\geq n/2.
 $$
 При
 $(n,q)\notin\{(4,2),(4,3)\}$ ложно неравенство $$n\geq \frac{q^t-1}{q-1}\geq  \frac{q^{n/2}-1}{q-1}.$$

 Отметим, что число $(q^t-1)/(q-1)$ одномерных подпространств $t$-мерного пространства $W$ над $\mathbb{F}_q$ как правило совпадает со степенью $\mu(\SL_t(q))$ минимального подстановочного представления группы $\SL(W)=\SL_t(q)$, см.~\cite{Cooperstein}.

  Для разбора случая $(-)$ мы детально рассмотрим действие $\SU(W)$ непосредственно на $\Omega\omega$ при $t=3$, чтобы исключить этот случай. При всех других $t$ мы воспользуемся информацией о степенях минимальных подстановочных представлений $\mu(\SU_t(q))$ группы $\SU_t(q)$, вычисленных в той же работе~\cite{Cooperstein}, а также цепочкой неравенств
  $$ 2t\geq n\geq m=|\Omega|\geq |\Omega\omega|\geq \mu(\SU_t(q)),$$
  которые имеют место аналогично случаю~$(1)$.

   Пусть $t=3$. По лемме Витта любые два изометричных подпространства пространства $W$ некоторый элемент $\SU(W)$ переводит друг в друга. Размерность вполне изотропных подпространств  в~$W$ не превосходит $[3/2]=1$. Вычисляя число соответствующих подпространств как индекс стабилизатора одного из них в группе $\SU(W)$, заключаем, что, вопреки неравенству $2t\geq |\Omega\omega|$, имеет место один из следующих случаев:
\begin{itemize}
  \item $V_i\omega$ вполне изотропны, их размерность  равна~1 и  $$|\Omega\omega|=q^3+1\geq 2^3
  +1= 9> 6=2t;$$
  \item $V_i\omega$ невырождены размерности~1,  $$|\Omega\omega|=q^2(q^2+q+1)\geq 2^2(2^2+2+1)= 28>6=2t;$$
  \item $V_i\omega$ невырождены размерности~2, их столько же, сколько (одномерных невырожденных) ортогональных дополнений к ним, и тем самым  имеет место оценка из предыдущего случая;
  \item $V_i\omega$ вырождены, но не изотропны, размерности~2, это, в точности, ортогональные дополнения к своим (одномерным изотропным) радикалам, их столько же,  сколько одномерных изотропных подпространств, т.е. справедлива оценка из первого случая.
\end{itemize}

Из \cite[таблица~1]{Cooperstein} с учетом изоморфизма $\SU_2(q)\cong \SL_2(q)$ имеем:
$$ \mu(\SU_t(q))= \left\{
\begin{array}{rl}
2, &\text{ если } (t,q)=(2,2), \\
3, &\text{ если }  (t,q)=(2,3), \\
5, &\text{ если }  (t,q)=(2,5), \\
7, &\text{ если }  (t,q)=(2,7), \\
6, &\text{ если }  (t,q)=(2,9), \\
11, &\text{ если }  (t,q)=(2,11), \\
q+1, &\text{ если }  t=2, q\ne 2,3,5,7,9, 11, \\
   2, &\text{ если }  (t,q)=(3,2), \\
  50, &\text{ если }  (t,q)=(3,5), \\
  q^3+1, &\text{ если }  t=3, q\ne 2,5, \\
  (q+1)(q^3+1), &\text{ если }  t=4, \\
  \displaystyle\frac{\left(q^{t-1}-(-1)^{t-1}\right)\left(q^{t}-(-1)^{t}\right)}{\left(q^{2}-1\right)}, &\text{ если }  t>4.
\end{array}
\right.
$$

Мы предполагаем, что $(t,q)\notin\{(2,2),(2,3)\}$, так как в противном случае $(n,q)\in\{(4,2),(4,3)\}$ и $m=n$, как в заключении утверждения~$(2)$.

Если $t=2$ и $q\in\{5,7,9,11\}$, то $\mu(\SU_t(q))>4=2t$, вопреки неравенству $2t\geq \mu(\SU_t(q))$. Если же $q\ne 5,7,9,11$, то снова
$$\mu(\SU_t(q))=q+1\geq 4+1>4=2t.$$

Случай $t=3$ исключен выше.

При $t=4$ также
$$
\mu(\SU_t(q))=(q+1)(q^3+1)\geq (2+1)(2^3+1)=27>8=2t.
$$

Наконец, исключим случай $t>4$. При четном $t=2d>4$ имеем
$$
\mu(\SU_t(q))=\displaystyle\frac{\left(q^{t-1}-(-1)^{t-1}\right)\left(q^{t}-(-1)^{t}\right)}{\left(q^{2}-1\right)}=(q^{t-1}+1)\sum\limits_{i=0}^{d-1}q^{2i}\geq d(2^{t-1}+1)>t\cdot 2^{t-2}>2t.
$$
При нечетном $t=2d+1>4$ имеем
$$
\mu(\SU_t(q))=\displaystyle\frac{\left(q^{t-1}-(-1)^{t-1}\right)\left(q^{t}-(-1)^{t}\right)}{\left(q^{2}-1\right)}=(q^{t}+1)\sum\limits_{i=0}^{d-1}q^{2i}\geq d(2^t+1)>(t-1)\cdot 2^{t-1}>2t.
$$
Всюду получаем противоречие с неравенством $2t\geq \mu(\SU_t(q))$. Утверждение~$(2)$ доказано.

%
%
%
%
%

Докажем~$(3)$. Пусть $x=x^{\phantom{x}}_Ux^{\phantom{x}}_W$, где $x^{\phantom{x}}_U\in \GL^\varepsilon(U)$ и $x^{\phantom{x}}_W\in \GL^\varepsilon(W)$. Так как элемент $x$ полупростой, имеет простой порядок и
  при $\varepsilon=+$ не имеет собственных векторов, a при $\varepsilon=-$ не имеет невырожденных собственных векторов в~$V$, имеем $$|x|=|x^{\phantom{x}}_U|=|x^{\phantom{x}}_W|\text{ и }(|x|, q-\varepsilon)=1.$$
 Мы не будем предполагать, что группа $L$ уже дана, но построим ее для элемента $x$ специальным образом, после чего убедимся, что она удовлетворяет условию леммы.

Из условия ясно, что $t\geq n/2\geq 2$. Так как $x$~--- полупростой элемент и $(n,q)\notin\{(4,2),(4,3),(4,5) \}$, по лемме~\ref{alpha_classic} существуют $m$ элементов $g_1,\dots,g_m\in \SL^\varepsilon(W)$ таких, что $$L^{\phantom{x}}_W:=\langle x_W^{g_1}, \dots, x_W^{g_m}\rangle\geq \SL^\varepsilon(W).$$ При этом $L^{\phantom{x}}_W=\langle \SL^\varepsilon(W), x^{\phantom{x}}_W\rangle$. Положим $$L^{\phantom{x}}_U:=\langle x^{\phantom{x}}_U\rangle\text{ и } L:=
\langle x^{g_1}, \dots, x^{g_m}\rangle.$$ Тогда $L\leq L^{\phantom{x}}_U\times L^{\phantom{x}}_W$ и проекции $L$ на сомножители $L^{\phantom{x}}_U$ и $L^{\phantom{x}}_W$ сюръективны. Чтобы показать, что $L$ удовлетворяет условию леммы, осталось установить, что $L$ совпадает с $L^{\phantom{x}}_U\times L^{\phantom{x}}_W$.

Для этого, во-первых, заметим, что $x^{\phantom{x}}_W\in \SL^\varepsilon(W)$; в частности, $L^{\phantom{x}}_W=\SL^\varepsilon(W)$. Это следует из того, что порядок образа элемента $x^{\phantom{x}}_W$ в факторгруппе $\GL^\varepsilon(W)/\SL^\varepsilon(W)$ делит $(|x|, q-\varepsilon)=1.$ Далее, рассмотрим канонический эпиморфизм
$$\phi: {L^{\phantom{x}}_U\times L^{\phantom{x}}_W}\rightarrow ({L^{\phantom{x}}_U\times L^{\phantom{x}}_W})\,/\,\Phi({L^{\phantom{x}}_U\times L^{\phantom{x}}_W}).$$ Учитывая, что $L^{\phantom{x}}_U$~--- группа простого порядка, $L^{\phantom{x}}_W\cong \SL^\varepsilon(W)$, и по известным свойствам подгруппы Фраттини \cite[гл.~A, лемма~(9.4)]{DH}, имеем $$\Phi({L^{\phantom{x}}_U\times L^{\phantom{x}}_W})=\Phi(L^{\phantom{x}}_U)\times \Phi(L^{\phantom{x}}_W)=\Phi(L^{\phantom{x}}_W)=\mathrm{Z}(L^{\phantom{x}}_W).$$ Поэтому $L^\phi$ изоморфна подгруппе в $\mathbb{Z}_{|x|}\times L_n^\varepsilon(q)$, проекции которой на каждый из сомножителей сюръективны. Но $L^\phi$ обладает композиционными факторами, изоморфными $\mathbb{Z}_{|x|}$ и $L_n^\varepsilon(q)$, поэтому $$L^\phi=({L_U\times L_W})\,/\,\Phi({L_U\times L_W}),\text{ откуда }L={L^{\phantom{x}}_U\times L^{\phantom{x}}_W}.$$

Теперь по доказанным утверждениям $(1)$ и $(2)$ леммы существует ${g\in\SL^\varepsilon(V)}$ такой, что $$G:=\langle L, x^{g}\rangle=\langle x^{g_1}, \dots, x^{g_m},x^{g}\rangle$$ неприводима и примитивна как подгруппа в~$\GL^\varepsilon(V)$. Так как $\SL^\varepsilon(W)\leq G$, группа $G$ содержит длинную корневую подгруппу группы $\SL^\varepsilon(V)$. Пусть $R$~--- нормальное замыкание в $G$ этой корневой подгруппы. Ввиду того, что подгруппа $G$ примитивна, подгруппа $R$ неприводима. В частности, $\Oo_p(R)=1$. Теперь как вытекает из \cite[теорема~II]{Kantor}, группа $R$ изоморфна одной из групп, указанных в утверждении~$(3)$ (для случая, когда ${\varepsilon=+}$, можно воспользоваться также результатами работ~\cite{McL1,McL2}).
\end{proof}

\begin{Lemma}\label{Field_Aut} {\rm \cite[предложение~4.9.1]{GLS3}} {\it Пусть $L={}^d\Sigma(q)$ --- простая группа лиева типа над полем ${\mathbb F}_q$, где $\Sigma$~--- неприводимая корневая система, $d$~--- либо пустой символ, либо $2$ (т.е. ${}^d\Sigma(q)\ne{}^3D_4(q)$). Пусть $x$  и $y$ --- автоморфизмы группы $L$, имеющие один и тот же простой порядок,  и допустим, что $x$ и $y$ являются одновременно полевыми или графово-полевыми автоморфизмами по модулю группы $\L$. Тогда подгруппы $\langle x\rangle$ и
$\langle y\rangle$ сопряжены элементом из~$\L$. Если $x$~--- полевой автоморфизм, то $${}^d\Sigma(q^{1/|x|})\leq C_{L}(x)\leq C_{\L}(x)=\widehat{{}^d\Sigma(q^{1/|x|})}.$$ Если же $x$~---
графово-полевой автоморфизм и ${}^d\Sigma\in\{A_{n-1},D_{n}\}$, то $|x|=2$ и $C_{L}(x)={}^2\Sigma(q^{1/2})$.}
\end{Lemma}

\begin{Lemma}\label{GraphAutGLU} {\rm \cite[лемма~1.7]{YRV}}
{\it Пусть $L=L_n^\varepsilon(q)$~--- простая проективная специальная линейная или унитарная группа и $n\geq 4$, и мы рассматриваем графовые по модулю $\L$ инволюции в группе автоморфизмов группы~$L$. Тогда справедливы следующие утверждения:
\begin{itemize}
\item[{\em (1)}] Если $n$ нечётно, то все графовые по модулю $\L$ инволюции сопряжены относительно $\L$, и каждая такая инволюция нормализует в $L$ некоторую подгруппу~$K$, изоморфную $\SL^\varepsilon_{n-2}(q)$, и индуцирует нетривиальный автоморфизм на~$K/\mathrm{Z}(K)$.
\item[{\em (2)}] Если $n$ чётно, а $q$ нечётно, то  существует три класса $\L$-сопряженности графовых по модулю $\L$ инволюций, и для их представителей  $x_0,x_+,x_-$ выполнено
    $$
    \mathrm{F}^*\left(C_L(x_\delta)\right)\cong\left\{
    \begin{array}{cc}
      S_n(q), &\text{если } \delta=0, \\
      O^+_n(q), &\text{если } \delta=+, \\
      O^-_n(q), &\text{если } \delta=-.
    \end{array}
    \right.
    $$
     Элемент $x_\delta$ нормализует подгруппу $K_{\delta}$ в $L$, которая 
     является образом относительно естественного гомоморфизма подгруппы вида $$\left(\GL_m^\varepsilon(q)\times \GL_{n-m}^\varepsilon(q)\right)\cap \SL_n^\varepsilon(q)$$
     из $\SL^\varepsilon_n(q)$, где
     $$
     m=\left\{\begin{array}{cl}
                2 & \text{ при } (\varepsilon,\delta)\in \{(+,0),(+,+),(-,+),(-,-)\} \\
                1 & \text{ при } (\varepsilon,\delta)\in \{(+,-),(-,0)\},
              \end{array}
     \right.
     $$
         и $x_\delta$ индуцирует нетривиальный автоморфизм на той компоненте факторгруппы $K_{\delta}/\mathrm{Z}(K_{\delta})$, которая изоморфна $\mathrm{PSL}_{n-m}^\varepsilon(q)$.
\item[{\em (3)}] Если оба числа $n$ и $q$ чётны, то  существуют два класса $\L$-сопряжённости графовых по модулю $\L$ инволюций. Для любой  инволюции $x$  в группе $L$ имеется $x$-инвариантная подгруппа $K$, для которой $K/\mathrm{Z}(K)\cong\PSL_n^\varepsilon(q)$, и $x$ индуцирует нетривиальный автоморфизм на~$K/\mathrm{Z}(K)$.
\end{itemize}}
\end{Lemma}

\begin{proof}
Информацию о классах $\L$-сопряженности инволюций, содержащихся в $\L\tau$, их числе, их представителях и централизаторах можно найти в \cite[лемма~10]{GuestLevy}.

Если $q$ нечетно, то все инволюции в $\L\tau$ сопряжены с $\tau$ посредством элементов из $\L$. Мы можем считать, что $x=\tau$.
Тогда $x$ нормализует подгруппу $K$ в $L$, которая является образом в $\mathrm{PSL}_n^\varepsilon(q)$ подгруппы вида $$\left(\GL_{n-1}^\varepsilon(q)\times \GL_1^\varepsilon(q)\right)\cap \SL_n^\varepsilon(q)$$  и состоит из образов матриц вида
$$
\left(
\begin{array}{cc}
A&0\\ 0&(\det A)^{-1}
\end{array}
\right),
$$
где $A$ пробегает $ \GL^\varepsilon_{n-1}(q)$, и нормализует подгруппу $K^{\infty}\cong \SL_{n-1}^\varepsilon(q)$. Ясно, что $x=\tau$ индуцирует нетривиальный автоморфизм на $K^{\infty}/Z\left(K^{\infty}\right)\cong \mathrm{PSL}^\varepsilon_{n-1}(q).$ Утверждение $(1)$ доказано.

Рассмотрим случай, когда $n$ чётно, а $q$ нечётно. При $\varepsilon=-$ доказательство утверждения о существовании требуемой подгруппы $K_\delta$ см. в~\cite[стр.~288]{LieSaxl} и~\cite[стр.~43]{Liebeck}. Поэтому мы считаем, что $\varepsilon=+$. В этом случае (см. \cite[стр.~285]{LieSaxl}) элемент $x$ сопряжён посредством элемента из $\L$ с одной из трех попарно несопряженных инволюций $x_0,x_+,x_-$, индуцированных на $L$ элементами $$J^0\tau,\, J^+\tau,\, J^-\tau\in \langle\GL_n(q),\tau\rangle$$ соответственно, где
\begin{itemize}
  \item $J^0$~--- блочно-диагональная матрица, на диагонали которой стоят блоки $$\left(\begin{array}{rr}0&-1\\1&0 \end{array}\right),$$
  \item $J^+$~--- блочно-диагональная матрица, на диагонали которой стоят блоки $$\left(\begin{array}{rr}0&1\\1&0 \end{array}\right),$$
  \item $J^-$~--- блочно-диагональная матрица, на диагонали которой стоят $n/2-1$ блоков $$\left(\begin{array}{rr}0&1\\1&0 \end{array}\right)$$ и один блок $$\left(\begin{array}{rr}\mu&0\\0&1 \end{array}\right),$$ где $\mu\in \F_q$ таково, что $-
      \mu/2$ не является квадратом в~$\F_q$.
\end{itemize}
При этом централизаторы элементов $x_\delta$ такие, как указано в утверждении~$(2)$. Покажем, что $x_\delta$ нормализует некоторую подгруппу $K_\delta$ такую, как указано в~$(2)$.

Элемент $J^-\tau$ нормализует подгруппу  вида $$\left(\GL_{n-1}(q)\times \GL_1(q)\right)\cap \SL_n(q)$$ группы $\SL_n(q)$, поэтому $x_-$ нормализует подгруппу $K_-$, которая определена также, как подгруппа $K$ доказательстве утверждения~$(1)$, и понятно, что $x_-$ индуцирует на $K_-$ неединичный автоморфизм. Для случая, когда $x$ сопряжен с $x_-$ утверждение~$(2)$ доказано.

Элементы $J^0\tau$ и $J^+\tau$  нормализуют подгруппу в $\GL_n(q)$, состоящую из блочно-диагональных матриц вида
$$
\left(
\begin{array}{cc}
A&0\\ 0&B
\end{array}
\right),
$$
где $A\in \GL_{n-2}(q)$, $B\in\GL_2(q)$. Пересекая ее с $\SL_n(q)$ и рассматривая образ этого пересечения в $\mathrm{PSL}_n(q)$, получим группу $K_0=K_+$ инвариантную относительно $x_0$ и $x_+$ и такую, как указано в утверждении $(2)$, c $m=2$. Утверждение~$(2)$ доказано.

Наконец, рассмотрим случай, когда $q$ и $n$ четны. Здесь, в соответствии с \cite[замечание~11]{GuestLevy} и \cite[лемма~3.7]{Liebeck},  при $\varepsilon=-$ мы будем рассматривать $\GL^\varepsilon_n(q)=\GU_n(q)$ как группу матриц относительно некоторого упорядоченного базиса $e_1, \dots, e_m, f_m, \dots, f_1$ тех линейных преобразований векторного пространства размерности  $n =2m$ над полем $\F_{q^2}$, которые сохраняют эрмитову форму $(\,\cdot\,,\,\cdot\,)$, определенную равенствами $$(e_i, e_j) =(f_i, f_j) =0 \text{ и } (e_i, f_j) =\delta_{ij} \text{ для всех } i,j=1,\dots,m.$$

Тогда (см. \cite[лемма~10]{GuestLevy}) у группы $L$ имеется ровно два класса $\L$-сопряженности графовых по модулю $\L$ автоморфизмов порядка~$2$, и некоторые их представители индуцированы элементами
\begin{itemize}
  \item $\tau$ и $ J^0t\tau$ из $\langle\GL_n(q),\tau\rangle$ при $\varepsilon=+$, где матрица $J^0$ определена выше, а 
  $$
  t=\left(
  \begin{array}{ccccc}
    1 & 1 & & & \\
     & 1 & & & \\
     &  & 1 & & \\
     &  &   &\ddots & \\
     & & & & 1
  \end{array}
  \right);
  $$
  \item $\varphi_q$ и $t_0\varphi_q$ из $\langle\GU_n(q),\varphi_q\rangle$, при $\varepsilon=-$, где
  $$
  t_0=\left(
  \begin{array}{ccccc}
    1 &  & & & 1\\
     & 1 & & & \\
     &  & \ddots & & \\
     &  &   & 1& \\
     & & & & 1
  \end{array}
  \right).
  $$
  (корневой элемент максимальной высоты).
  \end{itemize}

%

  Пусть вначале $\varepsilon=+$. Тогда $\tau$ и $J^0t\tau$ нормализуют  в $\SL_n(q)$ подгруппу всех блочно-диагональных матриц вида
  $$\left(
  \begin{array}{cc}
    A &  \\
     & B
  \end{array}
  \right),\quad \text{ где } \quad A\in \GL_2(q),\quad B\in\GL_{n-2}(q) \quad \text{ и }\quad \det A\cdot\det B=1.
  $$
  Образ этой подгруппы в $\PSL_n(q)$ содержит нормальную подгруппу ${K\cong \SL_{n-2}(q)}$, инвариантную относительно автоморфизмов, индуцированных $\tau$ и $J^0t\tau$, причем оба автоморфизма нетривиальны на $K/\mathrm{Z}(K)$.

Пусть теперь $\varepsilon=-$. Элемент $t_0$ в группе $\GL_n(q^2)$ централизует $\langle\tau,\varphi_q\rangle$-ин\-вари\-ант\-ную подгруппу $G^*$, состоящую из матриц вида
$$\left(
  \begin{array}{ccc}
    1 & & \\
    & A &  \\
    & & 1
  \end{array}
  \right),\quad \text{ где } \quad A\in \SL_{n-2}(q^2)
  $$
и централизует в ней каждую подгруппу. В ней подгруппа $K^*=C_{G^*}(\tau\varphi_q)\cong \SU_{n-2}(q)$ инвариантна относительно $\varphi_q$ и $t_0$ и содержится в $\SU_n(q)$. Элементы  $\varphi_q$  и $t_0\varphi_q$ нормализуют в $\SU_n(q)$ подгруппу $$N_{\SU_n(q)}(K^*)\cong (\GU_{n-2}(q)\times\GU_2(q))\cap\SU_n(q).$$ Обозначим  через $K$ образ в группе $\PSL_n(q^2)$ подгруппы $K^*$. Ясно, что $K\leq \PSU_n(q)$. Элементы  $\varphi_q$  и $t_0\varphi_q$ индуцируют на изоморфных группах $K^*/\mathrm{Z}(K^*)$ и $K/\mathrm{Z}(K)\cong \PSU_{n-2}(q)$ нетривиальные автоморфизмы с согласованным действием.
\end{proof}

Ввиду того, что в случае графового по модулю $\L$ автоморфизма группы $L=L_n^\pm(q)$ при $n=4$ оценка на $\alpha(x,L)$ является нерегулярной (см. лемму~\ref{alpha_classic}), нам при $n=4$ понадобится дополнительная информация о подгруппах в $L$ нормализуемых, но не централизуемых~$x$. Эту информацию дает следующая лемма, в которой мы пользуемся известными изоморфизмами $L_4^\pm(q)\cong O_6^\pm(q)$, см.~\cite[предложение~2.9.1]{KL}.

\begin{Lemma}\label{Graph_Inv_D_n} {\it Пусть $V$~--- векторное пространство над полем $\mathbb{F}_q$ нечетного порядка,  $\dim V=6$ и $V$ снабжено невырожденной симметрической билинейной формой знака $\varepsilon\in\{+,-\}$. Пусть ${\rm O}$~--- группа изометрий пространства $V$,  $\Delta$~--- его группа подобий, ${\rm SO}$~--- группа элементов из ${\rm O}$ c определителем $1$,  $\Omega={\rm O}'=\Delta'$. Пусть $$\overline{\phantom{x}}:\Delta\rightarrow\Delta/\mathrm{Z}(\Delta)$$
обозначает канонический эпиморфизм. Пусть $L=\overline{\Omega}=O_{2n}^\varepsilon(q)$. Тогда справедливы следующие утверждения.

\begin{itemize}
  \item[$(1)$] Канонический графовый автоморфизм $\overline{\gamma}$ группы $L$ содержится в $\overline{{\rm O}}\setminus \overline{{\rm SO}}$, а группа $\overline{\Delta}$ совпадает с $\langle\L,\overline{\gamma}\rangle$.
  \item[$(2)$] Все графовые по модулю $\L$ инволюции являются образами инволюций из $\Delta$.

  \item[$(3)$] Имеется три класса $\L$-сопряженности графовых по модулю $\L$ инволюций с представителями $\overline{\gamma}_1=\overline{\gamma}, \overline{\gamma}_2$ и $\overline{\gamma}_{2}'$, где $\gamma_i$ для каждого $i=1,2$~--- инволюция в ${\rm O}$, у которой собственное значение $-1$ имеет кратность ${2i-1}$, а $\gamma_{2}'$--- инволюция в $\Delta\setminus {\rm O}$, у которой кратность собственного значения $-1$ равна $3$. При этом каждая инволюция $\overline{\gamma}_1$ и $\overline{\gamma}_1$ нормализует, но не централизует подгруппу в $L$, изоморфную $O_{5}(q)\cong S_4(q)$, а инволюция $\overline{\gamma}_{2}'$ нормализует, но не централизует подгруппу в $L$, изоморфную $O_{4}^-(q)\cong L_2(q^2)$.

\end{itemize}
}
\end{Lemma}
\begin{proof} Все утверждения леммы, кроме существования подгрупп, нормализуемых, но не централизуемых соответствующими  инволюциями, следуют\footnote{В таблицах 4.5.1 и 4.5.2 информацию о графовых инволюциях и их централизаторах нужно брать в строке для групп $A^\varepsilon_m(q)$ при нечетном $m$, имея ввиду, что $A^\varepsilon_m(q)=D^\varepsilon_m(q)\cong O_6^\varepsilon(q)$.} из \cite[теоремы 4.5.1, 4.5.2, таблицы  4.5.1, 4.5.2 и замечание 4.5.4]{GLS3}.

Докажем, что инволюции $\gamma_i$, $i=1,2$ нескалярно действуют на некотором невырожденном подпространстве $U$ пространства $V$ таком, что $\dim U=5$. Поскольку $\gamma_i$~--- изометрия, подпространства $V_+$ и $V_-$, состоящие из собственных векторов с собственными значениями $1$ и $-1$ соответственно, ортогональны друг другу и $V= V_+\oplus V_-$. Поэтому  $V_+$ и $V_-$~--- невырожденные $\gamma_i$-инвариантные подпространства, на каждом из которых $\gamma_i$ действует скалярно. Пусть  $u\in V_+$~--- невырожденный вектор. Тогда подпространство $W=u^\perp$ является $\gamma_i$-инвариантным невырожденным дополнением к $\langle u\rangle$. Поскольку $\dim V_+=6-2i+1>1$, преобразование $\gamma_i$ имеет на $W$ оба собственных значения $1$ и~$-1$ и действует на $W$ нескалярно. Значит, $\overline{\gamma}_i$ нормализует, но не централизует образ коммутанта группы изометрий пространства~$W$, изоморфный~$O_{5}(q)$.

Согласно~\cite[таблица~4.5.1]{GLS3} централизатор $\overline{\gamma}_{2}'$ в $\L$ изоморфен некоторой группе  автоморфизмов группы  $O_{4}^-(q)\cong L_2(q^2)\cong O_3(q^2)$, содержащей эту группу, имеет тривиальный центр и нечетный индекс в~$\L$. Из леммы~\ref{guest} заключаем, что $\overline{\gamma}_{2}'$ нормализует, но не централизует подгруппу в $L$, изоморфную $O_{4}^-(q)$.
  \end{proof}

\begin{Lemma} \label{r_not divi} {\rm\cite[лемма~1.18]{YRV}}
Пусть  $L$~--- простая классическая группа или знакопеременная группа, и  $r$~--- нечётное простое число,
не делящее порядок группы $L$.
Тогда справедливы следующие утверждения.
\begin{itemize}
  \item[$(1)$] Если $L=A_n$, то $r\geq n+1$ и $n\leq r-1$.
  \item[$(2)$] Если $L=L_n(q)$, то  $ r\geq n+2$  и $n\leq r-2$.
  \item[$(3)$] Если $L=U_n(q)$, то  $ r\geq n+2$  и $n\leq r-2$.
  \item[$(4)$] Если $L=S_{n}(q)$, то  $r\geq n+3$   и $n\leq r-3$.
  \item[$(5)$] Если $S=O_{n}(q)$, $n$ нечётно, то $r\geq n+2$ и $n>r-2$.
  \item[$(6)$] Если $S=O^+_{n}(q)$ или $S=O^-_{n}(q)$, $n$ чётно, то  $r\geq n+1$ и $r\geq n+2$ соответственно.
  \item[$(7)$] Если $L=L_n(q^2)$, то     $r\geq 2n+3.$

\end{itemize}

\begin{Lemma}\label{Guest1}
Пусть $L=L_n^\varepsilon(q)$ и $x$~--- автоморфизм простого порядка группы $L$. Тогда имеет место одно из следующих утверждений:
\begin{itemize}
  \item[$(1)$] Существует элемент $g\in L$ такой, что подгруппа $\langle x,x^g\rangle$ неразрешима.
  \item[$(2)$] $q=3$, $x$~--- трансвекция и существуют $g_1,g_2\in L$ такие, что подгруппа $\langle x,x^{g_1},x^{g_2}\rangle$ неразрешима.
  \item[$(3)$] $\varepsilon=-$, $q=2$, $x$~--- псевдоотражение порядка $3$ и существуют $g_1,g_2,g_3\in L$ такие, что подгруппа $\langle x,x^{g_1},x^{g_2},x^{g_3}\rangle$ неразрешима.
\end{itemize}
\end{Lemma}

\begin{proof}
См.~\cite[теорема~A*]{Gu}.
\end{proof}

\begin{Lemma} \label{beta_A_n_prop} {\rm\cite[предложение~2]{YRV}}
Пусть  $L=A_n$, $n\geq 5$, $r\leq n$~--- нечетное простое число, и $x\in \Aut (L)$~--- элемент простого порядка. Тогда
\begin{itemize}
\item[$(1)$] $\beta_{r}(x,L)=r-1$, если $x$~--- транспозиция;
\item[$(2)$] справедливо одно из следующих утверждений:
 \begin{itemize}
\item[{\rm (a)}] $\beta_{r}(x,L)\leq r-1$;
\item[{\rm (б)}] $r=3$, $n=6$, $x$~--- инволюция, не лежащая в $S_6$, и   $\beta_{r,L}(x)= 3$.
\end{itemize}\end{itemize}
\end{Lemma}

\end{Lemma}

\section{Доказательство основной теоремы}

Доказательство основной теоремы  разобьём на несколько случаев. Начнем с разбора случаев малой размерности, которые рассматриваются в серии  последовательных лемм. Затем рассмотрим внутренне-диагональные автоморфизмы $x$, которые допускают естественное действие на ассоциированном векторном пространстве. Мы рассмотрим отдельно ситуации, когда $x$ стабилизирует одномерное подпространство т.\,н. унаследованного типа (в частности, здесь будет полностью разобран случай, когда автоморфизм $x$ индуцирован унипотентным элементом), и когда $x$~--- полупростой элемент, не имеющий таких одномерных инвариантных подпространств. И, наконец, разберем ситуации, когда автоморфизм $x$ является полевым, графово-полевым или графовым по модулю группы внутренне-диагональных автоморфизмов.

\subsection{Случаи малой размерности}

\begin{Lemma}\label{Small4}
Теорема {\em \ref{main}} справедлива, если $L\in\{L^\pm_4(2),L^\pm_4(3),L^\pm_4(5)\}$.
\end{Lemma}

\begin{proof}
Теорема~\ref{main} верна для группы $L_4(2)\cong A_8$ в силу леммы~\ref{beta_A_n_prop}, поэтому считаем, что $(\varepsilon,q)\ne(+,2)$. Рассмотрим случай $r=3$. Тогда по лемме~\ref{Guest1} за исключением случая, когда $L=U_4(2)$, а $x$~--- псевдоотражение порядка $3$, некоторые два или три сопряженных с $x$ элемента порождают неразрешимую группу. Ни одна из групп в условии леммы не имеет секций, изоморфных простым группам Сузуки, поэтому из теоремы Томпсона--Глаубермана \cite[гл.~II, следствие~7.3]{Glauberman} следует, что порядок такой неразрешимой подгруппы делится на 3. Таким образом, во всех случаях $\beta_3(x,L)\leq 3$.
Кроме того, по лемме~\ref{alpha_classic} имеем $\alpha(x,L)\leq 4$, за исключением следующих случаев:
\begin{itemize}
  \item $L=L_4(q)$, $q\in\{3,5\}$, $x$~--- графовый автоморфизм по модулю $\L$ и $\alpha(x,L)\leq 6$;
  \item $L=U_4(2)$, $x$~--- трансвекция и $\alpha(x,L)\leq 5$;
  \item $L=U_4(q)$, $q\in\{2,3,5\}$, $x$~--- графовый автоморфизм по модулю $\L$ и $\alpha(x,L)\leq 6$.
  \end{itemize}
  Отсюда вытекает справедливость леммы при $r>5$, а также, кроме перечисленных исключительных случаев, при $r=5$.
  Рассмотрим оставшиеся исключения.

  Пусть $L=L_4^\pm(5)$ и $x$~--- графовая инволюция по модулю $\L$. Тогда согласно лемме~\ref{GraphAutGLU}~(2) элемент $x$ нормализует некоторую подгруппу $K$ группы $L$ такую, что $K/\mathrm{Z}(K)\cong L_3^\pm(5)$ или $K/\mathrm{Z}(K)\cong L_2^\pm(5)$ и $x$ индуцирует на $K/\mathrm{Z}(K)$ нетривиальный автоморфизм. Отсюда и из лемм~\ref{MainL2q}, \ref{MainL3q} и~\ref{MainU3} заключаем, что $\beta_5(x,L)\leq 4$.

  Пусть $x$~--- графовый автоморфизм группы $L_4(3)$ по модулю внут\-рен\-не-диа\-го\-наль\-ных автоморфизмов и $|x|=2$. В таблице характеров группы $\Aut(L_4(3))$ из \cite{atlas} этому случаю соответствуют инволюции классов сопряженности $2D$, $2E$, $2F$ и $2G$, причем классы $2D$ и $2F$ переставляются диагональным автоморфизмом, поэтому $2E$ можно не рассматривать. Воспользуемся известным фактом из теории характеров, утверждающим, что для данных элементов $a,b$ и $c$ группы $G$ число $\mathrm{m}(a,b,c)$ пар $(u,v)\in a^G\times b^G$ таких, что $uv=c$, может быть найдено из таблицы характеров по формуле
  $$
  \mathrm{m}(a,b,c)=\frac{|a^G||b^G|}{|G|}\sum\limits_{\chi\in\Irr(G)}\frac{\chi(a)\chi(b)\overline{\chi(c)}}{\chi(1)},
  $$
  см. \cite[упр.~(3.9) на стр.~45]{Isaacs}.
  Используя таблицу характеров из~\cite{atlas} соответствующего расширения $L_4(3)$,  убеждаемся\footnote{Все вычисления производились с помощью системы GAP~\cite{GAP}.}, что в каждом из классов  $2F$ и $2G$ существует пара элементов, произведение которых принадлежит классу $5A$ и имеет порядок $5$ (более точно, $\mathrm{m}(2F,2F,5A)=\mathrm{m}(2G,2G,5A)=20$). Поэтому $\beta_5(x,L)=2$, если $x$~--- инволюция из классов $2F$ и $2G$. Также в классе $2A$, состоящем из внутренних инволюций есть пара элементов, произведение которых принадлежит $5A$, так как $\mathrm{m}(2A,2A,5A)=5$. При этом каждый элемент из $2A$ является произведением двух инволюций из $2D$, поскольку $\mathrm{m}(2D,2D,2A)=2$. Следовательно, для инволюции $x$ из $2D\cup 2E$ имеем $\beta_5(x,L)\leq 4=5-1$.

  Рассмотрим оставшиеся случаи.
  \begin{itemize}
    \item $L=U_4(2)$, $x$~--- трансвекция (класс $2A$):
    \begin{itemize}
      \item[] $\mathrm{m}(2A,2A,2B)=2, \,\, \mathrm{m}(2B,2B,5A)=5  \,\,\,\Rightarrow \,\,\,\beta_5(x, U_4(2))\leq 4, \text{ если }  x\in 2A.$
    \end{itemize}
    \item $L=U_4(2)$, $x$~--- графовая по модулю $\L$ инволюция (классы $2C$, $2D$):
    \begin{itemize}
    \item[] $ \mathrm{m}(2C,2C,5A)=5 \quad \Rightarrow\quad \beta_5(x, U_4(2))= 2, \text{ если }  x\in 2C;$
      \item[] $\mathrm{m}(2D,2D,2B)=2, \,\, \mathrm{m}(2B,2B,5A)=5  \,\,\, \Rightarrow \,\,\, \beta_5(x, U_4(2))\leq 4, \text{ если }  x\in 2D.$
    \end{itemize}
    \item $L=U_4(3)$, $x$~--- графовая по модулю $\L$ инволюция (классы $2D$, $2E$, $2F$):
    \begin{itemize}
      \item[] $\mathrm{m}(2D,2D,2A)=2, \, \mathrm{m}(2A,2A,5A)=5 \,\,\, \Rightarrow\,\,\, \beta_5(x, U_4(3))\leq 4, \text{ если }  x\in 2D;$
      \item[] $ \mathrm{m}(2E,2E,5A)=5 \quad \Rightarrow\quad \beta_5(x, U_4(3))= 2, \text{ если }  x\in 2E;$
      \item[] $ \mathrm{m}(2F,2F,5A)=5 \quad \Rightarrow\quad \beta_5(x, U_4(3))= 2, \text{ если }  x\in 2F.$
    \end{itemize}

  \end{itemize}
\end{proof}

\begin{Lemma}\label{MainL2q}
Теорема {\em \ref{main}} справедлива, если $L= L_2(q)$.
\end{Lemma}

\begin{proof} Ведем индукцию по $q$. Справедливость теоремы \ref{main} в группах $L_2(4)$, $L_2(5)$ и $L_2(9)$ следует из леммы~\ref{beta_A_n_prop} и изоморфизмов $L_2(4)\cong L_2(5)\cong A_5$ и $L_2(9)\cong A_6$. Поэтому считаем, что
 $q\not=4,5,9$.

 Если $x$ не является полевым автоморфизмом порядка $2$, то в силу  леммы~\ref{alpha_classic} для любого простого делителя $s$ порядка группы~$S$ имеем $$\beta_s(x)\leq\alpha(x)\leq 3,$$ откуда следует утверждение леммы. В частности, лемма верна, если $q$~--- простое число.

Предположим теперь, что $\vert x\vert=2$ и автоморфизм $x$ полевой по модулю $\widehat{L}$. 
Рассмотрим $\varphi=\varphi_{{q^{1/2}_{\phantom{0}}}}$~--- канонический полевой автоморфизм порядка 2, индуцированный отображением
$$\left(\begin{array}{cc}
u&v\\
w&z
\end{array}\right)\mapsto
\left(\begin{array}{cc}
u^{q^{1/2}_{\phantom{0}}}&v^{q^{1/2}_{\phantom{0}}}\\
w^{q^{1/2}_{\phantom{0}}}&z^{q^{1/2}_{\phantom{0}}}
\end{array}\right).
$$
Поскольку в $\langle \widehat{L},\varphi\rangle$ подгруппы $\langle x\rangle$ и $\langle \varphi\rangle$ сопряжены  относительно $\widehat{L}$ по лемме~\ref{Field_Aut}, можно считать, что $x=\varphi$. Так как $q\not=9$, по лемме~\ref{alpha_classic} имеем $\beta(x,L)\leq\alpha(x,L)\leq 4$. Поэтому для любого $r\geq 5$ справедливо утверждение теоремы~\ref{main}.

Предположим,  $r=3$. Тогда $r$ делит порядок $L$ и ${s=r=3}$. Кроме того,  по лемме~\ref{Field_Aut} $$\mathrm{PSL}_2\big({q^{1/2}_{\phantom{0}}}\big)\leq C_L(x)\leq \mathrm{PGL}_2\big({q^{1/2}_{\phantom{0}}}\big)$$ и $C_L(x)$~---  подгруппа четного индекса в~$L$. Она почти проста, так как $q\ne 4,9$, и следовательно $\mathrm{Z}(C_L(x))=1$. По лемме~\ref{guest} элемент $x$ нормализует, но не централизует некоторую подгруппу, сопряженную с $\Oo^{p'}(C_L(x))\cong L_2\big({q^{1/2}_{\phantom{0}}}\big)$. Остается воспользоваться предположением индукции.
%
\end{proof}

\begin{Lemma}\label{MainU3}
Теорема {\em \ref{main}} справедлива, если $L= U_3(q)$.
\end{Lemma}

\begin{proof}
Из леммы \ref{alpha_classic} следует, что $\alpha(x)\leq 3$ за исключением случая, когда $q=3$ и  $\vert x\vert=2$ и $\alpha(x)=4$. Поэтому $$\beta_s(x,L)\leq \alpha(x)\leq 4,$$ откуда следует лемма для $r>3$, и, более того,
 $$\beta_s(x,L)\leq \alpha(x)\leq 3,$$ кроме случая, когда $r=3$,  $L=U_3(q)$ и $|x|=2$. Для разбора последнего случая (в котором с необходимостью $s=r=3$) воспользуемся таблицей характеров группы $U_3(3)$ и ее группы автоморфизмов из \cite{atlas}. В $U_3(3)$ все инволюции образуют один класс сопряженности $2A$, а все не внутренние~--- класс~$2B$.
Рассуждая, как в заключительной части леммы~\ref{Small4} и используя GAP \cite{GAP}, находим:
\begin{itemize}
      \item[] $ \mathrm{m}(2A,2A,3B)=3\,\, \quad \Rightarrow\quad \beta_3(y, U_3(3))= 2, \text{ если }  y\in 2A;$
      \item[] $ \mathrm{m}(2B,2B,3A)=36 \quad \Rightarrow\quad \beta_3(y, U_4(3))= 2, \text{ если }  y\in 2B.$
    \end{itemize}
\end{proof}

\begin{Lemma}\label{MainL3q}
Теорема {\em \ref{main}} справедлива, если $L=L_3(q)$.
\end{Lemma}

\begin{proof}
 Пусть вначале автоморфизм $x$ не является графово-полевым по модулю $\widehat{L}$. В силу леммы~\ref{alpha_classic} имеем $\alpha(x)\leq 3$. Значит,  для любого простого делителя $s$ порядка группы $ L$ справедливо неравенство $\beta_s(x,L)\leq \alpha(x,L)\leq 3$, откуда следует утверждение леммы.

Предположим теперь, что $x$ графово-полевой автоморфизм по модулю $\widehat{L}$ и $|x|=2$, т.е. $\widehat{L}x=\widehat{L} \tau\varphi$, где $$\tau:A\mapsto (A^{-1})^\top,\quad\text{ a }$$
$$\varphi=\varphi_{q^{1/2}}:\left(a_{ij}\right)\mapsto \left(a_{ij}^{q^{1/2}}\right).$$  По лемме~\ref{Field_Aut} можно считать, что $x=\tau\varphi$. В этом случае  $\alpha(x)\leq 4$ согласно лемме~\ref{alpha_classic}, т.\,е. утверждение леммы справедливо для любого $r\geq 5$. Предположим, что $r=3$.  Тогда, $s=r=3$ и, как следует из леммы~\ref{Field_Aut}, $$C_{\widehat{L}}(x)\cong \mathrm{PGU}_3(q^{1/2})\text{ и }\Oo^{p'}(C_L(x))\cong U_3(q^{1/2}).$$ Поскольку 
$$\Oo^{p'}(C_L(x))\leq C_{{L}}(x)\leq C_{\widehat{L}}(x),$$ легко убедиться, что индекс $\vert {L:C_L(x)}\vert$ чётен. Кроме того, $\mathrm{Z}(C_L(x))=1$. В силу леммы~\ref{guest}, существует подгруппа $M$ группы $L$ такая, что $M$ сопряжена с $C_L(x)$ и $x$ нормализует, но не централизует $M$. Пусть $y$~--- автоморфизм группы ${O^{p'}(M)\cong U_3(q^{1/2})}$, индуцированный~$x$. Остается воспользоваться леммой~\ref{MainU3} и тем, что $$\beta_3(x,S)\leq \beta_3(y,\Oo^{p'}(M))\leq 3.$$
%
%
\end{proof}

\begin{Lemma}\label{beta5S4q}
 { Если $L=S_4(q)$, где $q$ нечетно, и  
 $x\in \L$~--- инволюция, то
$
\beta_5(x,L)\leq 4.
$}
\end{Lemma}

\begin{proof} Воспользуемся изоморфизмами
$$L\cong O_5(q)\cong \Omega_5(q),\quad \L \cong \mathrm{SO}_5(q)\quad\text{и}\quad \mathrm{GO}_5(q)=\langle -E\rangle\times \mathrm{SO}_5(q),
$$
где $E$~--- единичная матрица $5\times 5$, см. \cite[таблица~2]{atlas}, \cite[предложение 2.9.1]{KL}. Возьмем $\eta\in\{1,-1\}$ так, чтобы выполнялось сравнение $$q\equiv \eta\pmod 4.$$
Рассмотрим 5-мерное векторное пространство $V$ над $\F_q$ с определенной на нем невырожденной симметрической билинейной  формой $(\,\cdot\,,\,\cdot\,)$ и некоторым базисом $e_1,\dots,e_5$ таким, что для некоторого $\zeta\in\F_q^*$ выполнены равенства
$$
(e_i,e_j)=\zeta\delta_{ij}\text{ для } i,j\in\{1,\dots,5\}.
$$
Как известно, $\zeta\in(\F_q^*)^2$, если $\eta=1$ и $\zeta\notin(\F_q^*)^2$, если $\eta=-1$, и знак формы равен~$\eta$.

Пусть $G$~--- группа всех линейных преобразований $g$ пространства $V$ таких, что
$$
(ug,vg)=(u,v)\quad\text{для всех}\quad u,v\in V
$$
и отождествим $\L$ с подгруппой $\{g\in G\mid \det g=1\}$ индекса 2 в группе $G$. Тогда $L=G'$.
Выберем в $\langle e_1,e_5\rangle$ ортогональные друг другу ненулевые векторы $e_1',e_5'$ так, чтобы $$(e_1',e_1')=(e_5',e_5')\notin (\F_q^*)^2\zeta. $$ Тогда для соответствующих ${i\ne j}$ сужение формы на $\langle e_i,e_j\rangle$ имеет знак $\eta$, а на $\langle e_i,e_j'\rangle$~--- знак $-\eta$. Обозначим для краткости полученные упорядоченные базисы следующим образом:
$$\begin{array}{cc}
    \mathcal{E}: & e_1,e_2,e_3,e_4,e_5, \\
    \mathcal{E}': & e_1',e_2,e_3,e_4,e_5'.
  \end{array}
$$
Матрицы линейного преобразования $g$ пространства $V$, записанные в базисах $\mathcal{E}$ и $\mathcal{E}'$ соответственно, будем обозначать символами
$$
[g]\quad\text{и}\quad [g]'.
$$

Известно \cite[таблица~4.5.1]{GLS3}, что в группе $L$ имеется два класса сопряженности инволюций и еще два класса имеется в $\L\setminus L$. Укажем их представители. Это $x_1^\square$,
$x_1^\boxtimes$, $x_2^\square$ и
$x_2^\boxtimes$ такие, что
$$\begin{array}{c}
    \left[x_1^\square\right]=\left[x_1^\boxtimes\right]'=\mathrm{diag}\,(-1,-1,\phantom{-}1,\phantom{-}1,\phantom{-}1), \\
    \left[x_2^\square\right]=\left[x_2^\boxtimes\right]'=\mathrm{diag}\,(-1,-1,-1,-1,\phantom{-}1).
  \end{array}
$$
Таким образом определенные элементы легко раскладываются в произведение отражений относительно векторов $e_i$ и $e_i'$, легко определяются их спинорные нормы, откуда заключаем, что
$$
x_1^\square,x_2^\square\in L,\quad \text{a} \quad x_1^\boxtimes,x_2^\boxtimes\in\L\setminus L.
$$

Из определения следует, что элементы $x_1^\boxtimes$ и $x_2^\boxtimes$ стабилизируют невырожденное подпространство $$W=\langle e_2,e_3,e_4,e_5'\rangle,$$
ограничение формы на которое имеет знак ``$-$''. Кроме того, эти элементы индуцируют на $W$ нескалярные преобразования $y_1$ и $y_2$, сохраняющие форму на $W$. Поэтому если $H$~--- группа изометрий пространства $W$, то $y_1$ и $y_2$ индуцируют неединичные автоморфизмы на $$H'/\mathrm{Z}(H')\cong O_4^-(q)\cong L_2(q^2),$$
порядок группы $L_2(q^2)$ делится на 5, и по лемме~\ref{MainL2q} имеем
$$
\beta_5(x_i^\boxtimes, L)\leq \beta_5(y_i, H'/\mathrm{Z}(H'))\leq 4\quad\text{для}\quad i=1,2.
$$

Далее, рассмотрим вложение $S_5\hookrightarrow G$, задаваемое действием группы $S_5$ на индексах векторов $e_1,\dots,e_5$. Образ подстановки $\sigma\in S_5$ обозначим через $\sigma^*$. Инволюция $\tau=(15)(24)\in A_5$ инвертирует в $A_5$ элемент $\sigma=(12345)$ порядка $5$, поэтому $\tau\tau^\sigma=\sigma^2$~--- элемент порядка $5$ и
$\beta_5(\tau, A_5)=2$. Из простоты группы $A_5$ следует, что $\sigma^*,\tau^*\in L$, и $\beta_5(\tau^*,L)=2$. Так как $\tau^*\in L$, инволюция $\tau^*$ сопряжена с $x_1^\square$ или~$x_2^\square$, а из того, что кратность собственного значения $-1$ у преобразования $\tau^*$ равна $2$, заключаем, что $\tau^*$ сопряжена с~$x_1^\square$. Таким образом,
$$
\beta_5(x_1^\square,L)=2.
$$
Наконец, для элементов $g=(15)(34)^*$ и $h=(25)(34)^*$ из $L$ имеем $$(x_2^\square)^g(x_2^\square)^h=x_1^\square,$$ откуда
$$
\beta_5(x_2^\square,L)\leq 2\beta_5(x_1^\square,L)= 4
$$ по лемме~\ref{estim}.
\end{proof}

\subsection{Доказательство теоремы~\ref{main}: структура и общие замечания}\label{general}


Мы начинаем {\it Доказательство теоремы}~\ref{main}. Пусть, как в условии теоремы, $L=L^\varepsilon_n(q)$ и $x\in\Aut(L)$~--- автоморфизм простого порядка. Напомним также, что $r$~--- нечетное простое число, и $s\in\pi(L)$ выбрано так, что $s=r$, если $r\in\pi(L)$ и $s>r$, если $r\notin\pi(L)$. 
Наша цель~--- доказать неравенство
\begin{equation}\label{ineq}
 \beta_s(x,L)\leq \left\{
                            \begin{array}{lr}
                              r, & \text{если } r=3, \\
                              r-1, & \text{если } r>3.
                            \end{array}
                          \right.
\end{equation}
Ведем рассуждения индукцией по~$|L|$.

В силу лемм~\ref{Small4}--\ref{MainL3q} мы можем считать, что $n\geq 4$ и при $n=4$ число    $q$ отлично от $2,3,5$. Можно считать также, что за исключением ситуации, когда $n=4$ и $x$~--- графовый по модулю $\L$ автоморфизм, выполнено неравенство
\begin{equation}\label{rn-inequality}
r\leq n,
\end{equation}
 поскольку в противном случае по лемме~\ref{alpha_classic} имеем
$$
\beta_s(x, L)\leq\alpha(x,L)\leq n\leq r-1,
$$ и неравенство~(\ref{ineq}) верно. Отсюда и из леммы~\ref{r_not divi}  вытекает, что
\begin{equation}\label{s=r}
r\text{ делит }|L|\text{ и, следовательно, }s=r.
\end{equation}

В разделах \ref{indiag(1)} и \ref{indiag(varepsilon)} мы рассмотрим все случаи, когда $x\in\L$, а в \ref{field-graph-and}~--- все случаи, когда $x\in\Aut(L)\setminus \L$.

Поскольку $x$ имеет простой порядок, в случае $x\in\L$ элемент $x$ либо унипотентен, либо полупрост.  Отождествим группу  $L^+_n(q)=L_n(q)$ с $\mathrm{PSL}(V)$, где $V$~--- $n$-мерное векторное пространство  над полем $\F_{q}$, а группу $L^-_n(q)=U_n(q)$ с $\mathrm{PSU}(V)$, где $V$~--- $n$-мерное векторное пространство c невырожденной эрмитовой формой над $\F_{q^2}$. Тогда $\L=\mathrm{PGL}^\varepsilon(V).$ Рассматривая элемент $x\in\mathrm{PGL}^\varepsilon(V)$, мы будем говорить, что $x$-инвариантное подпространство $U\leq V$ {\it имеет унаследованный тип}, если $U$~--- произвольное подпространство при $\varepsilon=+$ или $x$ унипотентном, и $U$ невырождено при $\varepsilon=-$ и $x$ полупростом. Разделы \ref{indiag(1)} и~\ref{indiag(varepsilon)} соответствуют случаям, когда для $x$ существует и не существует одномерное подпространство унаследованного типа.

Будем также говорить, что $x\in\mathrm{PGL}^\varepsilon(V)$ {\it действует скалярно} на $x$-инвариантном подпространстве~$U$, если у некоторого прообраза в $\GL^\varepsilon(V)$ элемента $x$ сужение на $U$ пропорционально тождественному преобразованию   (эквивалентно, любое одномерное подпространство в $U$ инвариантно относительно $x$).

\subsection{Внутренне-диагональный автоморфизм,\\ стабилизирующий одномерное подпространство\\ унаследованного типа}\label{indiag(1)}

В этом разделе доказательства мы разберем все случаи, когда $x$ оставляет инвариантным некоторое одномерное подпространство $U$ унаследованного типа пространства $V$ (в частности, такое предположение охватывает случай, когда элемент $x$ унипотентен).  
Возможны следующие подслучаи:
 \begin{itemize}
   \item[(а)] $x$ унипотентен и $\varepsilon=+$;
   \item[(б)] $x$ унипотентен и $\varepsilon=-$;
   \item[(в)] $x$ полупрост и либо $\varepsilon=+$, либо $U$ невырождено.
 \end{itemize}
 Рассмотрим случай (а). Пусть $P_1$~--- стабилизатор некоторой прямой в естественном модуле, а $P_2$~--- стабилизатор гиперплоскости. По лемме~\ref{Parabolic}  $x\in P_i\setminus \Oo_p(P_i)$ для $i=1$ или $i=2$ с точностью до сопряжения элементом из $L$, причём элемент $x$ на компонентах фактора $P_i/\Oo_p(P_i)$ действует нескалярно.  В частности, это означает, что $x\notin \Oo_\infty(P_i)$. Рассмотрим канонический эпиморфизм $$\overline{\phantom{G}}:P_i\rightarrow P_i/\Oo_\infty(P_i).$$ Тогда $\overline{x\vphantom{P_1}}\ne1$ и $\mathrm{F}^*(\overline{P_i})\cong L_{n-1}(q)$. Из леммы~\ref{r_not divi}  и неравенства~(\ref{rn-inequality}) вытекает, что $r$ делит  $|L_{n-1}(q)|$, и по  лемме \ref{inductionLemma}
 $$\beta_r(x, L)\leq\beta_r\big(\overline{\vphantom{P_1}x}, \mathrm{F}^*(\overline{P_i})\big).$$ Применив предположение индукции, получаем  неравенство~(\ref{ineq}). Для случая (а) теорема доказана.

Пусть имеет место случай (б).

Если  $x$ стабилизирует также некоторое невырожденное одномерное подпространство пространства~$V$, то справедливость теоремы~\ref{main} устанавливается повторением рассуждений для случая (а) с естественной заменой $L^+_{n-1}(q)=L_{n-1}(q)$ на $L^-_{n-1}(q)=U_{n-1}(q)$ и параболической подгруппы $P_i$ на стабилизатор невырожденного одномерного подпространства, изоморфный образу в $\mathrm{PGU}_n(q)$ подгруппы $$\GU_1(q)\times \GU_{n-1}(q)\leq \GU_n(q).$$ Поэтому в случае (б) считаем, что $x$ не стабилизирует никакое одномерное подпространство.

 Рассмотрим две параболические максимальные подгруппы
\begin{itemize}
\item $P_1$~--- стабилизатор максимального вполне изотропного подпространства $W$ пространства~$V$ (по~лемме Витта размерность такого подпространства равна~$[n/2]$) и
  \item $P_2$~--- стабилизатор изотропного одномерного подпространства $U$ из~$V$.
  \end{itemize}
  По лемме~\ref{Parabolic} для $i=1$ или $i=2$ с точностью до сопряжения элементом из $L$ имеем  $x\in P_i\setminus \Oo_p(P_i)$. Как и в случае (а), имеем $x\notin \Oo_\infty(P_i)$ и рассматриваем канонический эпиморфизм $$\overline{\phantom{G}}:P_i\rightarrow P_i/\Oo_\infty(P_i).$$

Допустим, $i=1$, т.\,е. $x$ стабилизирует максимальное вполне изотропное подпространство $W$ и индуцирует на нем нетождественное преобразование.
Тогда $$\mathrm{F}^*\left(\overline{P_1}\right)\cong L_{[n/2]}(q^2),\quad\text{причем}\quad \overline{x\vphantom{P_1}}\in \mathrm{F}^*\left(\overline{P_1}\right)^\sharp.$$ Неравенство~(\ref{rn-inequality}) и малая теорема Ферма показывают, что $r$ делит $\left|L_{[n/2]}(q^2)\right|$, поскольку
$$
2\left[{n}/{2}\right]\geq n-1\geq r-1.
$$
 Теперь лемма  \ref{inductionLemma} в сочетании с предположением индукции дает:
 $$
 \beta_r(x,L)\leq
   \beta_r\big(\overline{x\vphantom{P_1}},\mathrm{F}^*\big(\overline{P_1}\big)\big)\leq \left\{
                            \begin{array}{lr}
                              r, & \text{если } r=3, \\
                              r-1, & \text{если } r>3,
                            \end{array}
                          \right.$$

  Допустим, $i=2$. Тогда $\mathrm{F}^*\big(\overline{P_2}\big)\cong U_{n-2}(q)$, причем $\overline{x}\in \mathrm{F}^*\big(\overline{P_1}\big)^\sharp$.
  Если $r$ делит $|U_{n-2}(q)|$, то снова в соответствии с леммой~\ref{inductionLemma} имеем
  $$
 \beta_r(x,L)\leq
   \beta_r\big(\overline{x\vphantom{P_1}},\mathrm{F}^*\big(\overline{P_2}\big)\big),$$ откуда по предположению индукции получаем неравенство~(\ref{ineq}).
  Поэтому, разбирая оставшиеся подслучаи случая (б), считаем, что число  $r$, не делит $|U_{n-2}(q)|$, а это по лемме~\ref{r_not divi} влечет неравенство $r\geq n$.  Поэтому  $n=r$ в силу неравенства~(\ref{rn-inequality}). В частности, $n$ нечетно и $n\geq 5$.

В соответствии с леммой~\ref{UnipotentUnitary} имеет место один из следующих подслучаев
\begin{itemize}
  \item $x$  стабилизирует невырожденное подпространство размерности~$1$; этот подслучай уже исключен;
  \item $x$ стабилизирует максимальное вполне изотропное подпространство $W$ и индуцирует на нем нетождественное преобразование; этот подслучай также исключен;
  \item $q$ нечетно и существует сопряженный с $x$ элемент $x^g$ такой, что в подгруппе $\langle x, x^g\rangle$ некоторая инволюция\footnote{Строго говоря, применяя лемму~\ref{UnipotentUnitary}, мы должны говорить о прообразах элементов $x$ и $x^g$ в группе $\SU_n(q),$ и инволюция также берется в $\SU_n(q)$, однако нескалярность действия инволюции на максимальном вполне изотропном подпространстве влечет нескалярность ее действия на всем $V$, а значит ее образ в $U_n(q)$ по-прежнему будет инволюцией.} $y$ стабилизирует максимальное вполне изотропное подпространство $W$, индуцируя на нем нескалярное преобразование; именно этот подслучай нам осталось исключить для завершения разбора случая~(б).
\end{itemize}
Как и выше, пусть $P_1$~--- стабилизатор максимального вполне изотропного подпространства $W$, и мы рассматриваем канонический эпиморфизм $$\overline{\phantom{G}}:P_1\rightarrow P_1/\Oo_\infty(P_1).$$ При этом $\mathrm{F}^*(\overline{P_1})\cong L_{[n/2]}(q^2)=L_{(r-1)/2}(q^2)$. Инволюция $\overline{y}$ нормализует, но не централизует эту подгруппу, индуцируя на $L_{(r-1)/2}(q^2)$ внутренне-диагональный автоморфизм. Из малой теоремы Ферма следует, что $r$ делит $q^{r-1}-1=(q^2)^{(r-1)/2}-1$, а значит, делит $|\mathrm{F}^*(\overline{P_1})|$. В частности, $s=r$.

 Если $n=r=5$, то $\dim W=(r-1)/2=2$. По лемме~\ref{betaL2q} в этом случае
 $$
 \beta_s(y,L)=\beta_r(y,L)\leq\beta_r\left(\overline{y\vphantom{P_1}}, \mathrm{F}^*(\overline{P_1})\right)=2,
 $$
откуда, учитывая, что $y\in\langle x,x^g\rangle$, по лемме~\ref{estim} получаем $$\beta_s(x,L)\leq 2\beta_s(y,L)\leq 4=r-1,$$ т.\,е. доказываемая теорема верна.

Пусть $n=r>5$. Тогда $r\geq 7$,  $\dim W=(r-1)/2\geq 3$ и по лемме~\ref{alpha_classic} с учетом того, что $\overline{y}$~--- внутренне-диагональная инволюция группы $L_{(r-1)/2}(q^2)$, имеем
$$
 \beta_r(y,L)\leq\beta_r\big(\overline{y\vphantom{P_1}}, \mathrm{F}^*(\overline{P_1})\big)\leq \alpha\big(\overline{y\vphantom{P_1}}, \mathrm{F}^*(\overline{P_1})\big)\leq \frac{r-1}{2},
 $$
откуда $$\beta_r(x,L)\leq 2\beta_s(y,L)\leq r-1.$$ Случай (б) разобран полностью.

Рассмотрим случай (в), когда полупростой элемент $x$ стабилизирует одномерное подпространство $U$, которое невырождено в случае $\varepsilon=-$. Пусть $W$~--- $x$-ин\-вари\-ант\-ное дополнение к $U$, причем $W\perp U$, если $\varepsilon=-$.

Не уменьшая общности рассуждений, мы можем считать, что прообраз элемента $x$ в  $\GL^\varepsilon(V)$ действует нескалярно на $W$. В самом деле, если $x$ действует скалярно на $W$, то $W$ состоит из собственных векторов прообраза $x^*$ элемента $x$ в $\GL^\varepsilon(V)$, причем, поскольку $x\ne 1$, собственное значение~$\lambda$, которому соответствуют вектора из $W$, отлично от собственного значения~$\mu$, которому соответствуют вектора из $U$. Выберем в $W$ одномерное подпространство $U_0$ (невырожденное, если $\varepsilon=-$) и рассмотрим $x$-инвариантное (ортогональное при $\varepsilon=-$) дополнение $W_0$ к нему. Так как $n\geq 4$, среди собственных значений ограничения $x^*$ на $W_0$ присутствуют как $\lambda$, так и $\mu$ и $x^*$ действует на $W_0$ нескалярно, и мы можем заменить $U$ на~$U_0$, а $W$ на~$W_0$.

Теперь элемент $x$ содержится в образе в $\mathrm{PGL}^\varepsilon_n(q)$ подгруппы вида $$\GL^\varepsilon_1(q)\times \GL_{n-1}^\varepsilon(q)$$ и индуцирует неединичный автоморфизм $\overline{x}$ на единственном неабелевом композиционном факторе $L^\varepsilon_{n-1}(q)$ этого образа. Как и при рассмотрении случая (а), на основании леммы~\ref{r_not divi} и неравенства~(\ref{rn-inequality}) заключаем, что
 $r$ делит $|L^\varepsilon_{n-1}(q)|$, откуда по  лемме~\ref{inductionLemma} и предположению индукции выводим неравенство~(\ref{ineq}). Для случая (в) теорема доказана.

\subsection{Внутренне-диагональный автоморфизм,\\ не стабилизирующий одномерных подпространств\\ унаследованного типа}\label{indiag(varepsilon)}

Здесь мы рассмотрим все случаи, когда полупростой элемент $x\in\L=\mathrm{PGL}^\varepsilon(V)$ не имеет  инвариантных одномерных  подпространств унаследованного типа.

Если $x$ действует неприводимо на $V$, то $$\beta_r(x,L)\leq\alpha(x,L)\leq 3\leq \left\{
                            \begin{array}{lr}
                              r, & \text{если } r=3, \\
                              r-1, & \text{если } r>3.
                            \end{array}
                          \right.$$
в силу леммы~\ref{alpha_irreducible}.
Поэтому считаем, что $V$ обладает собственным ненулевым $x$-ин\-ва\-риант\-ным подпространством. Сведем ситуацию к случаю, когда это подпространство унаследованного типа и, следовательно, размерности $\geq 2$, чтобы затем воспользоваться леммой~\ref{irredlinear}. По лемме~\ref{cases} если у $x$ нет собственных инвариантных подпространств унаследованного типа, то $\varepsilon=-$, $n$ четно и $x$ стабилизирует вполне изотропное подпространство $U$ размерности $n/2$. Так как элемент $x$ полупрост, по теореме Машке $x$ стабилизирует также некоторое вполне изотропное подпространство $W$ той же размерности и такое, что
$$
V=U\oplus W.
$$
Если $x$ действует скалярно на обоих подпространствах $U$ и $W$, возьмем ненулевые векторы $u\in U$ и $w\in W$ такие, что $(u,w)\ne 0$. Для прообраза $x^*\in \GL^-(V)$ элемента $x$ подпространства $U$ и $W$ являются пространствами собственных векторов, причем, поскольку $x^*\notin \mathrm{Z}(\GL^-(V))$ векторы $u$ и $w$ отвечают разным собственным значениям. Тем самым, $\langle u,w\rangle$~--- собственное $x$-инвариантное подпространство унаследованного типа.

Предположим, что $x$ действует нескалярно на $U$ или $W$. Тогда $x^*$ содержится в стабилизаторе в $\GL^-(V)$ этого подпространства, изоморфном $\GL_{n/2}(q^2)$, элемент $x$ содержится в образе $H$ этого стабилизатора в $\L=\mathrm{PGU}(V)$ и индуцирует неединичный автоморфизм $\bar{x}$ на $H^\infty/\mathrm{Z}(H^\infty)\cong L_{n/2}(q^2)$. Далее, как обычно, из неравенства~(\ref{rn-inequality}) и леммы~\ref{r_not divi} заключаем, что $r$ делит $|L_{n/2}(q^2)|$ и  неравенство~(\ref{ineq}) верно по предположению индукции и в силу неравенства
$$\beta_r(x,L)\leq\beta_r(\bar{x},H^\infty/\mathrm{Z}(H^\infty)).$$

Если же $x$ действует скалярно как на $U$, так и на $W$, то на   $\langle u,w\rangle$ элемент $x$ действует нескалярно, как мы и хотели.

Итак, до конца данного раздела считаем, что $x\in \L$~--- полупростой элемент, $U$~--- ненулевое $x$-инвариантное подпространство унаследованного типа наименьшей размерности (в частности, $x$ действует неприводимо на $U$) и $W$~--- $x$-инвариантное подпространство также унаследованного типа, дополняющее $U$ до~$V$.  При этом ${\dim U\geq 2}$ и $t=\dim W=n-\dim U\geq \dim U$. Элемент $x$ действует нескалярно на $W$ и поэтому индуцирует автоморфизм $\bar{x}$ группы $\mathrm{PSL}^\varepsilon(W)\cong L^\varepsilon_t(q)$. Если $r$ делит $|L^\varepsilon_t(q)|$, то по предположению индукции из неравенства
$$\beta_r(x,L)\leq\beta_r(\bar{x},L^\varepsilon_t(q))$$
следует справедливость неравенства~(\ref{ineq}). Поэтому считаем, что $r$ не делит   $|L^\varepsilon_t(q)|$, в частности, $r\geq 5$. По лемме~\ref{r_not divi}     выполнено неравенство
$$
t\leq r-2.
$$
Из леммы~\ref{irredlinear} следует, что если
$$m=\left\{\begin{array}{cc}
                                                                                                                          t & \text{ при } t>2, \\
                                                                                                                          3 &   \text{ при } t=2,
                                                                                                                        \end{array}
  \right.$$
  то некоторые $m+1$ элементов, сопряженных с $x$ посредством элементов из $L$ порождают подгруппу $H$ в~$L$, в которой содержится нормальная подгруппа из следующего списка
\begin{itemize}
  \item $L^\varepsilon_n(q)$ при любых $q$,
  \item $S_n(q)$ при любых $q$,
  \item $L_n^-(q_0)$ при $q=q_0^2$,
  \item $O^\pm_n(q)$ при четных $q$,
  \item $\Sym_{n+1}$ при $q=2$ и нечетном $n> 6$,
  \item $\Sym_{n+2}$ при $q=2$ и четном $n\geq 6$.
\end{itemize}
Порядок подгруппы $H$ делится на $r$, как следует из леммы~\ref{r_not divi}. Следовательно,
$$
\beta_r(x,L)\leq m+1=\left\{\begin{array}{rc}
t+1\leq r-1 & \text{ при } t>2, \\
4=5-1\leq r-1 &   \text{ при } t=2.
                                                                                                                        \end{array}
  \right.
$$
Рассматриваемый случай разобран полностью.

\subsection{Полевые, графово-полевые и графовые автоморфизмы}\label{field-graph-and}

В группе $\SL_n^\varepsilon(q)$ рассмотрим подгруппу $H$, состоящую из матриц вида
$$
\left(
\begin{array}{cc}
  A &  \\
   & 1
\end{array}
\right),
$$
где $A$ пробегает группу $\SL_{n-1}^\varepsilon(q)$, и образ $K$ подгруппы $H$ в $L_n^\varepsilon(q)$.Ясно, что , что $K/\mathrm{Z}(K)\cong L_{n-1}^\varepsilon(q)$. Из неравенства~(\ref{rn-inequality}) и леммы~\ref{r_not divi} следует, что $r$ делит $|L_{n-1}^\varepsilon(q)|$. Ясно также, что $H$ и $K$ инвариантны относительно автоморфизмов  $\varphi_{p^m}$ и $\tau$, причем $y$ индуцирует на $K/\mathrm{Z}(K)$ нетождественный автоморфизм $\bar{y}$ в каждом из следующих случаев:
\begin{itemize}
  \item $\varepsilon=+$ и $y$ совпадает с $\tau$, $\varphi_{q_0}$, если $q=q_0^t$ для простого числа $t$, или $\tau \varphi_{q_0}$, если $q=q_0^2$;
  \item $\varepsilon=-$ и $y$ совпадает с  $\varphi_{q_0}$, если $q^2=q_0^t$ для простого числа $t$.
\end{itemize}
  В силу лемм \ref{Field_Aut} и \ref{GraphAutGLU}  если автоморфизм $x$ является полевым или графово-полевым по модулю $\L$ или же $n$ нечетно и $x$ является графовым по модулю $\L$,  подгруппа $\langle x\rangle$ сопряжена относительно $\L$ с $\langle y\rangle$ для одного из упомянутых $y$. По предположению индукции из соотношений
  $$
  \beta_r(x,L)=\beta_r(y,L)\leq \beta(\bar{y}, K/\mathrm{Z}(K))
  $$
  выводим для этих случаев неравенство~(\ref{ineq}).

Остается таким образом рассмотреть случай, когда $n$ четно и $x$ является графовым по модулю $\L$. При четных $n>4$ из  неравенства~(\ref{rn-inequality}) следует, что
$$
r\leq n-1.
$$ Это же неравенство верно, если $n=4$ и $r=3$.
 В этих случаях лемма~\ref{r_not divi} влечет, что порядки групп $L_{n-1}^\varepsilon(q)$ и даже $L_{n-2}^\varepsilon(q)$ делятся на~$r$. Из леммы~\ref{GraphAutGLU} вытекает, что любой графовый по модулю $\L$ автоморфизм $x$ группы $L$ нормализует, но не централизует подгруппу $H$ в $L$ такую, что
$H^\infty/\mathrm{Z}(H^\infty)$ изоморфна $L_{n-1}^\varepsilon(q)$ или $L_{n-2}^\varepsilon(q)$. Применяя предположение индукции к автоморфизму, индуцированному $x$ на $H^\infty/\mathrm{Z}(H^\infty)$, получаем справедливость неравенства~(\ref{ineq}).

Осталось рассмотреть случай, когда $n=4$ и $r>3$. Так как $\alpha(x,L)\leq 6$ по лемме~\ref{alpha_classic} при $r\geq 7$ имеем
$$\beta_s(x,L)\leq \alpha(x,L)\leq 6=7-1\leq r-1.$$ Поэтому считаем, что $r=5$, откуда $r$ делит~$|L|$ и $s=r$. Учитывая лемму~\ref{GraphAutGLU}, достаточно рассмотреть случай, когда $q$~--- простое число. Так как случаи $q=2,3,5$ уже рассмотрены в лемме~\ref{Small4}, считаем, что $q\geq 7$, в частности, $q$ нечетно.

 Воспользуемся изоморфизмом $$L=L^\varepsilon_4(q)\cong O^\varepsilon_6(q)$$ и рассмотрим $L$ как проективную ортогональную группу, а элемент $x$ как графововую по модулю группы внутренне-диагональных автоморфизмов инволюцию в этой группе.
 По лемме~\ref{Graph_Inv_D_n} инволюция $x$ нормализует, но не централизует в $L$ подгруппу $H$, изоморфную либо $O_5(q)\cong S_4(q)$, либо $O_3(q^2)\cong L_2(q^2)$, и в обоих случаях $|H|$ делится на~5.

 Если $H\cong S_4(q)$, то $\Aut(H)=\widehat{H}$, поэтому  $x$ индуцирует на $H$ внутренне-диа\-го\-наль\-ный автоморфизм $\bar{x}$ и по леммам~\ref{inductionLemma} и~\ref{beta5S4q} имеем
 $$
 \beta_5(x,L)\leq \beta_5(\bar{x},H)\leq 4=5-1.
 $$

 В случае $H\cong L_2(q^2)$ имеем те же неравенства со ссылкой на лемму~\ref{MainL2q}. \hfill{$\Box$}

 \section{Доказательство теоремы~\ref{Cor}}

Пусть собственное подмножество $\pi$ множества всех простых чисел содержит по крайней мере два различных элемента.  Пусть $r$~--- наименьшее простое число, не принадлежащее~$\pi$, и
$$
m=\left\{
\begin{array}{ll}
  r, & \text{если } r\in\{2,3\}, \\
  r-1, & \text{если } r>3.
\end{array}
\right.
$$
  Допустим, теорема~\ref{Cor} неверна. Тогда $r\geq 3$ по лемме~\ref{2notinpi}, и группа наменьшего порядка среди конечных групп, у которых неабелевы композиционные факторы изоморфны группам из множества
  \begin{multline*}
    \{A_n\mid n\geq 5\}\cup
    \{L^\varepsilon_n(q)\mid n\geq 2, \varepsilon=\pm , q\text{ --- степень простого числа}\\ \text{и } (\varepsilon, n,q)\ne(\pm,2,2), (\pm,2,3), (-,3,2) \},
  \end{multline*}
  согласно лемме~\ref{red} изоморфна группе $G$, содержащей нормальную подгруппу $L$ из указанного множества, и элемент простого порядка $x$ такие, что $L$ не является $\pi$- или $\pi'$-группой, $G=\langle L,x\rangle$ и любые $m$ сопряженных с $x$ элементов порождают $\pi$-подгруппу в $G$. Пусть $s$~--- наименьший простой делитель $|L|$, не принадлежащий $\pi$. Тогда либо $r$ делит $|L|$ и $s=r$, либо $r$ не делит $|L|$ и $s>r$. По теореме~\ref{main} и лемме~\ref{beta_A_n_prop} выполнено неравенство $$\beta_s(x,L)\leq m,$$ т.\,е. существуют элементы ${g_1,\dots,g_m\in L}$ такие, что $|\langle{ x^{g_1},\dots,x^{g_m}}\rangle|$ делится на $s$. Но это противоречит тому, что любые $m$ сопряженных с $x$ элементов порождают $\pi$-под\-группу.  \hfill{$\Box$}

\end{document}